\documentclass[12pt,reqno]{amsart} 
 
   \usepackage{amsthm} 
   \usepackage{amsmath} 
   \usepackage{amssymb} 
   \usepackage{amsfonts} 
   \usepackage{amscd} 
   \usepackage[mathscr]{eucal} 
   \usepackage{verbatim} 
 
\newtheorem{thm}{Theorem}[section] 
\newtheorem{lem}[thm]{Lemma} 
\newtheorem{prop}[thm]{Proposition} 
\newtheorem{cor}[thm]{Corollary} 
\newtheorem{rem}[thm]{Remark} 
\newtheorem{defi}[thm]{Definition}

\newcommand{\vlto}[1] 
{\,\begin{picture}(#1,3) 
\put(0,2){\vector(1,0){#1}} 
\end{picture}\,} 
 
\newcommand{\neto}[1]
{\,\begin{picture}(#1,3) 
\put(0,0){\vector(1,1){#1}} 
\end{picture}\,} 
 
\newcommand{\nto}[1]
{\,\begin{picture}(#1,3) 
\put(0,0){\vector(0,1){#1}} 
\end{picture}\,} 
 
\newcommand{\adots}{\mathinner{\mkern1mu\raise1pt\hbox{.} 
\mkern2mu\raise4pt\hbox{.} 
\mkern2mu\raise7pt\hbox{.}\mkern1mu}} 
 
\newcommand{\A}{{\bf A}} 
\newcommand{\B}{{\bf B}} 
\newcommand{\G}{{\bf G}} 
\newcommand{\wG}[1]{\ensuremath{{\mbox{\bf G}}}^\sim_{#1}}

\newcommand{\M}{{\bf M}} 
\newcommand{\N}{{\bf N}} 
\newcommand{\PP}{{\bf P}} 
 
\newcommand{\T}{{\bf T}} 
\newcommand{\U}{{\bf U}} 
\newcommand{\Z}{{\bf Z}}

\newcommand{\bN}{\overline{N}} 
\newcommand{\bP}{\overline{P}} 
 
\newcommand{\bk}{\overline{k}} 
\newcommand{\bn}{\overline{n}} 
 
\newcommand{\hG}{\widehat{G}} 
\newcommand{\hM}{\widehat{M}} 
 
\newcommand{\hP}{\widehat{P}} 
\newcommand{\hT}{\widehat{T}}

\newcommand{\cT}{\mathcal{T}} 
 
\newcommand{\Ad}{\mathbb A} 
\newcommand{\cx}{\mathbb C} 
 
\newcommand{\zl}{\mathbb Z} 
\newcommand{\rl}{\mathbb R}

\renewcommand{\aa}{\alpha} 
\renewcommand{\a}[1]{\ensuremath{{\alpha_{#1}}}} 
\newcommand{\ba}[1]{\ensuremath{{\overline{\alpha}_{#1}}}} 
\newcommand{\bb}{\beta} 
\renewcommand{\b}[1]{\ensuremath{{\beta_{#1}}}} 
\newcommand{\e}{\epsilon} 
\newcommand{\g}{\gamma} 
\newcommand{\gl}[1]{\ensuremath{\mbox{GL}_{#1}}} 
\newcommand{\gso}[1]{\ensuremath{\mbox{GSO}_{#1}}} 
\newcommand{\gsp}[1]{\ensuremath{\mbox{GSp}_{#1}}} 
\newcommand{\gspin}[1]{\ensuremath{\mbox{GSpin}_{#1}}} 
\newcommand{\wgspin}[1]{\ensuremath{\mbox{GSpin}}^\sim_{#1}} 
\renewcommand{\l}{\lambda} 
 
\newcommand{\oo}{\omega} 
 
\newcommand{\s}{\sigma} 
\renewcommand{\sl}[1]{\ensuremath{\mbox{SL}_{#1}}} 
\newcommand{\so}[1]{\ensuremath{\mbox{SO}_{#1}}} 
\renewcommand{\sp}[1]{\ensuremath{\mbox{Sp}_{#1}}} 
\newcommand{\spin}[1]{\ensuremath{\mbox{Spin}_{#1}}}

\newcommand{\w}[1]{\ensuremath{\widetilde{#1}}} 
\newcommand{\h}[1]{\ensuremath{\widehat{#1}}} 
 
\begin{document} 
 
\title{Generic Transfer for General Spin Groups} 
\author[Mahdi Asgari]{Mahdi Asgari$^\star$} 
\address{School of Mathematics \\ 
Institute for Advanced Study \\ 
Einstein Drive \\
Princeton, NJ 08540 \\
USA} 
\email{asgari@math.ias.edu} 
\thanks{$^\star$Partially supported by a Rackham Fellowship at the  
University of Michigan and the NSF grant DMS--0211133 at the Fields Institute.}  
	 
\author[Freydoon Shahidi]{Freydoon Shahidi$^\dag$} 
\address{Mathematics Department \\  
Purdue University \\ 
West Lafayette, IN 47907 \\ 
USA} 
\email{shahidi@math.purdue.edu} 
\thanks{$^\dag$Partially supported by the NSF grant DMS--0200325.}

\begin{abstract} 
We prove Langlands functoriality for the generic spectrum of 
general spin groups (both odd and even). Contrary to other recent  
instances of functoriality, our resulting automorphic  
representations on the general linear group will not be self-dual.  
Together with cases of classical groups, this completes the list  
of cases of split reductive groups whose $L$-groups have classical  
derived groups. The important transfer from $\gsp{4}$ to $\gl{4}$  
follows from our result as a special case.  
\end{abstract} 
 
\maketitle 
 
\section{Introduction}\label{intro} 
 
Let $\G$ be a connected reductive group over a number field $k$.  
Let $G = \G(\Ad)$, where $\Ad$ is the ring of adeles of $k$.  
Let $^LG$ denote the $L$-group of $G$ and fix an embedding  
\[ \iota : \ ^LG\hookrightarrow \gl{N}(\cx) \times W(\bk/k), \] 
where $W(\bk/k)$ is the Weil group of $k$. Without loss of  
generality we may assume $N$ to be minimal.  
Let $\pi=\otimes_v^\prime \pi_v$ be an automorphic representation  
of $G$. Then for almost all $v$, the local representation $\pi_v$  
is an unramified  
representation and its class is determined by a semisimple  
conjugacy class $[t_v]$ in $^LG$. Here $v$ is a finite place of $k$.  
Let $\Pi_v$ be the unramified representation of $\gl{N}(k_v)$  
determined by the conjugacy class $[\iota(t_v)]$ generated by $\iota(t_v)$.  
{\it Langlands' functoriality conjecture} then demands the existence of an  
automorphic representation $\Pi^\prime=\otimes_v^\prime \Pi_v^\prime$  
of $\gl{N}(\Ad)$ such that $\Pi^\prime_v \simeq \Pi_v$ for all  
the unramified places $v$. In this paper we prove functoriality in the cases  
where $\G$ is not classical but the derived group 
$^LG^0_D$ of the connected component  of its $L$-group is. 
(We follow the convention that a classical group is the  
stabilizer of a symplectic, orthogonal, or hermitian  
non-degenerate bilinear form. Hence, for example, spin  
groups would not be considered classical.) 
 
As we explain later, a major difficulty in establishing this 
result is the absence of any useful matrix representation when 
the groups themselves are not classical, the subject matter 
of the present paper, forcing us to use rather complicated 
abstract structure theory in order to prove stability of the 
corresponding root numbers. 

We shall be mainly concerned with quasi-split groups and those  
automorphic representations which are induced from generic cuspidal  
ones. The problem clearly reduces to establishing functoriality  
for generic cuspidal representations of $G=\G(\Ad)$.  
 
The cases when $\G$ is a quasi-split classical group were  
addressed in \cite{cogdell-kim-ps-sha,ckpss-classical,kim-muthu},  
unless $\G$ is a quasi-split special orthogonal group which should  
be taken up by the authors of \cite{ckpss-classical}.  
 
In this paper we will prove the functorial transfer of generic  
cuspidal representations when $\G=\gspin{m}$, the split general  
spin group of semisimple rank $n=[m/2]$. We describe their  
structure in detail in  
Section \ref{structure}. In particular, these are split reductive  
linear algebraic groups of type $B_n$ or $D_n$ whose derived groups  
are double coverings of split special orthogonal groups. Moreover,  
the connected component of their Langlands dual groups are  
$^LG^0=\gsp{2n}(\cx)$ or $\gso{2n}(\cx)$, respectively.  
Then, $^LG = \gso{2n}(\cx) \rtimes W(\bk/k)$ or  
$\gsp{2n}(\cx) \rtimes W(\bk/k)$ according to whether $m$ is  
even or odd. The map $\iota$ is the natural  
embedding. Observe that $^LG^0_D$ is now a classical  
group and these groups are precisely the ones for which $\G^0_D$ is  
not classical but $^LG^0_D$ is. The transfer is to the space of  
automorphic representations of $\gl{2n}(\Ad)$.  
 
It is predicted by the theory of (twisted) endoscopy of  
Kottwitz, Langlands, and Shelstad \cite{kot-shel,langlands-shel} that  
the representations of $\gl{2n}(\Ad)$ which are in the image  
of this transfer must be of the form  
\begin{equation}\label{omega}  
\Pi = \w{\Pi} \otimes \oo  
\end{equation} 
for some gr{\"o}ssencharacter $\oo$. If $\oo_\Pi$ is the  
central character of $\Pi$, this implies that $\oo_\Pi / \oo^n$  
must be a quadratic character $\mu$ of $k^\times \backslash \Ad^\times$.  
Each $\mu$ then determines a quadratic extension of $k$ via  
class field theory and the group $\G$ which has transfers to  
automorphic representations of the type just mentioned would be the  
quasi-split form $\gspin{2n}^*$ of $\gspin{2n}$ associated to the  
quadratic extension. The split case then corresponds to $\mu\equiv 1$  
which is the content of the present paper.

If a representation $\pi$ of $\gspin{2n}^*(\Ad)$ transfers to $\Pi$  
on $\gl{2n}(\Ad)$ satisfying $\Pi \simeq \w{\Pi}\otimes\oo$  
for some gr{\"o}ssencharacter $\oo$, then $\oo=\omega_\pi$ and  
$\oo_\Pi = \oo_\pi^n \mu$, where $\oo_\pi$ and $\oo_\Pi$ denote the  
central characters of $\pi$ and $\Pi$, respectively, and $\mu$ is a  
quadratic gr{\"o}ssencharacter associated with the quasi-split  
$\gspin{2n}^*$. While we are not able to show that every $\Pi$  
satisfying (\ref{omega}) is transfer of an automorphic representation  
$\pi$, we show that our transfers satisfy (\ref{omega}). (In fact,  
we will prove that $\Pi$ is {\it nearly equivalent} to $\w{\Pi}\otimes\oo$  
for now. See Theorem \ref{transfer}.) 
 
We should note here that if $\Pi$ is an automorphic transfer  
to $\gl{2n+1}(\Ad)$ satisfying (\ref{omega}), then $\oo = \theta^2$ for  
some $\theta$ and $\Pi\otimes\theta^{-1}$ is then self-dual. Therefore,  
this is already subsumed in the self-dual case which is a  
case of standard twisted endoscopy. On the other hand, the case of  
$\gl{2n}(\Ad)$ discussed above is an example of the most general  
form of transfer that twisted endoscopy can handle.  
 
As explained earlier, in this paper we prove Langlands' functoriality  
conjecture in the form discussed, for all generic cuspidal representations  
of split $\gspin{m}(\Ad)$. In other words, we establish generic transfer  
from $\gspin{m}(\Ad)$ to $\gl{2n}(\Ad)$. Extension of this transfer  
to the non-generic case would require  
either the use of models other than Whittaker models  
or of Arthur's twisted trace formula. As far as we know new models  
for these groups have not been developed and the fact that these  
groups are not classical may make the matters complicated. On the other  
hand, the use of Arthur's twisted trace formula is at present depending  
on the validity of the fundamental lemmas which are not available  
for these groups. We refer to \cite{arthur-gsp4} for information on the case of
$\gsp{4}$. 

To state our main theorem, fix a Borel subgroup $\B$ in $\G$ and a  
maximal (split) torus $\T$ in $\B$, and denote  
the unipotent radical of $\B$ by $\U$. Let $\psi$ be a non-trivial continuous  
additive character of $k\backslash\Ad$. As usual, we use $\psi$ and a  
fixed splitting (i.e., the choice of Borel above along with a collection  
of root vectors, one for each simple root of $\T$, cf. page 13 of \cite{kot-shel},  
for example) to define a non-degenerate additive character of  
$\U(k)\backslash\U(\Ad)$, again denoted by $\psi$. (Also see Section 2 of  
\cite{shahidi-imrn}).   
 
Let $(\pi,V_{\pi})$ be an irreducible cuspidal automorphic representation  
of $\G(\Ad)$. The representation $\pi$ is said to be {\it globally generic}  
if there exists a cusp form $\phi\in V_{\pi}$ such that  
\begin{equation}  
\int_{\U(k)\backslash\U(\Ad)} \phi(ng)\psi^{-1}(n) dn \not= 0.  
\end{equation} 
Note that cuspidal automorphic representations of general linear  
groups are always globally generic. Two irreducible automorphic  
representations $\Pi$ and $\Pi^\prime$ of $\gl{N}(\Ad)$ are said to be   
{\it nearly equivalent} if there is a finite set of places $T$ of $k$  
such that $\Pi_v\simeq\Pi^\prime_v$ for all $v\not\in T$. Our main result  
is the following:  
 
\begin{thm}\label{transfer} 
Let $k$ be a number field and let $\pi=\otimes^\prime \pi_v$ be an  
irreducible globally generic cuspidal automorphic representation of either  
$\gspin{2n+1}(\Ad)$ or $\gspin{2n}(\Ad)$. Let $S$ be a non-empty  
finite set of non-archimedean places $v$ such that for $v\not\in S$  
we have that $\pi_v$ and $\psi_v$ are unramified. Then, there exists an  
automorphic representation $\Pi=\otimes^\prime\Pi_v$ of $\gl{2n}(\Ad)$  
such that for all archimedean $v$ and all non-archimedean $v$ with  
$v\not\in S$  the homomorphism parameterizing the local representation  
$\Pi_v$ is given by  
\[ \Phi_v = \iota\circ\phi_v : W_v\longrightarrow\gl{2n}(\cx), \]  
where $W_v$ denotes the local Weil  
group of $k_v$ and $\phi_v : W_v\longrightarrow\,^LG^0$ is the  
homomorphism parameterizing $\pi_v$.  
Moreover, if $\oo_\Pi$ and $\oo_\pi$ denote  
the central characters of $\Pi$ and $\pi$, respectively, then  
$\oo_\Pi = \oo_\pi^n$. Furthermore, if $v$ is an archimedean  
place or a non-archimedean place with $v\not\in S$, then  
$\Pi_v \simeq \w{\Pi}_v \otimes \oo_{\pi_v}$. In particular,  
the representations $\Pi$ and $\w{\Pi}\otimes\omega_\pi$ are  
nearly equivalent. 
\end{thm} 
 
At the non-archimedean places $v$ where $\pi_v$  
is unramified with the semisimple conjugacy class $[t_v]$ in  
$^LG^0$ as its Frobenius-Hecke (or Satake) parameter, this amounts to  
the fact that the local representation $\Pi_v$ is the unramified  
irreducible admissible representation determined by the conjugacy  
class generated by $\iota(t_v)$ in $\gl{2n}(\cx)$.

Our method of proof is that of applying an appropriate version  
of converse theorem \cite{converse1,converse2} to a family of  
$L$-functions whose required properties, except for one, are  
proved in \cite{shahidi:90annals,gel-sha,kim-sha,kim-cjm}. The  
exception, i.e., the main stumbling block for applying the  
converse theorem, is that of stability of certain root numbers under  
highly ramified twists. In \cite{shahidi-imrn} the root numbers, or  
more precisely the local coefficients, were expressed as a Mellin transform  
of certain Bessel functions. Applying this to our case requires a good  
amount of development and calculations. This is particularly important  
since the necessary Bruhat decomposition of these groups are more  
complicated than the classical groups. For that we have to resort to the  
use of abstract theory of roots which are harder since no reasonable  
matrix representation is available for these groups. Moreover, the main  
theorem in \cite{shahidi-imrn} is based on certain assumptions whose  
verification requires our calculations.  
 
The fact that $\gspin{2n}$ has a disconnected center makes matters  
even more complicated. This led us to use an extended group $\wgspin{2n}$  
of $\gspin{2n}$ so that our proof of stability proceeds smoothly.

There are two important transfers that are special cases of this theorem.  
The first is the generic transfer from $\gsp{4}=\gspin{5}$  
to $\gl{4}$ whose proof, as far as we know, has never been published  
before. We should point out that even the unpublished proofs of this 
result  are based on methods that are fairly disjoint from ours. 
We finally remark that our result in this case also gives an 
immediate proof of the holomorphy of spinor L-functions for generic 
cusp forms on $\gsp{4}$ (cf. Remark \ref{ramin}).

The second sepcial case is when $\G=\gspin{6}$. In this case our transfer 
gives the exterior square transfer from  
$\gl{4}$ to $\gl{6}$ due to H. Kim \cite{kim-ex2} which, when composed  
with symmetric cube of a cuspidal representation on $\gl{2}(\Ad)$,  
leads to its symmetric fourth. 
 
The issue of strong transfer, which has been successfully treated in the  
cases of classical groups thanks to existence of descent from $\gl{n}$ to classical  
groups \cite{grs,sou}, still needs to wait until the descent or other  
techniques are established  
for representations of $\gl{2n}(\Ad)$ which satisfy (\ref{omega}). 
 
Further applications such as non-local estimates towards the Ramanujan  
conjecture as well as some of the other applications established in  
\cite{cogdell-kim-ps-sha,ckpss-classical}  
will be addressed in future papers. As mentioned earlier, the cases of  
quasi-split $\gspin{}$ groups is the subject matter of our next paper.  
 
Here is an outline of the contents of each section.  
In Section \ref{str} we review the structure theory of the groups  
involved in this paper. In particular, we give detailed description of 
the root data for $\gspin{}$ groups and their extensions. We then prove  
the necessary analytic  
properties of local $L$-functions in Section \ref{local}. In particular,  
we discuss standard module conjecture which is another local ingredient.  In  
Section \ref{stability} we go on to prove the most crucial local result,  
i.e., stability of $\g$-factors under twists by highly ramified characters.  
This is where we do the calculations with root data mentioned above and  
use the extended group.  
We then prove the necessary analytic properties of the global $L$-functions  
in Section \ref{global} which will prepare us to apply the converse theorem  
in Section \ref{proof}. In Section \ref{apps} we include the special cases  
mentioned above along with some other local and global consequences of  
the main theorem. 
 
The authors would like to thank J. Cogdell for helpful discussions. The  
first author would also like to thank G. Prasad for answering many of  
his questions on algebraic groups and M. Reeder for helpful discussions.  
Some of this work was done while the authors were visiting the Fields  
Institute as part of the Thematic Program on Automorphic Forms in Spring  
of 2003. We would like to thank the Fields Institute for their hospitality  
and support. 
 
\section{Structure Theory}\label{str} 
 
We review the structure theory for the families of algebraic groups relevant  
to the current work, namely, $\gspin{2n+1}$ and $\gspin{2n}$ as well as their  
duals $\gsp{2n}$ and $\gso{2n}$. We will also introduce the group $\wgspin{2n}$  
which is closely related to $\gspin{2n}$. It shares the same derived group as  
that of $\gspin{2n}$. However, contrary to $\gspin{2n}$ which has disconnected  
center, the center of $\wgspin{2n}$ is connected. We will need this group for  
our purposes as we will explain later.  
 
\subsection{Root data for $\gspin{}$ groups}\label{structure} 
 
We first describe the algebraic group $\gspin{}$ and its standard Levi  
subgroups in terms of their root data. We will heavily rely on these description in  
the computations of Section \ref{stability}. 
 
Let $\G=\gspin{m}$, where $m=2n+1$ or $m=2n$. We now describe the root datum for $\G$. 
 
\begin{defi} Let  
\[ X={\mathbb Z}e_0\oplus{\mathbb Z}e_1\oplus\dots\oplus{\mathbb Z}e_n \]  
and  
\[ X^\vee={\mathbb Z}e_0^*\oplus{\mathbb Z}e_1^*\oplus\dots\oplus{\mathbb Z}e_n^* \]  
and let $\langle \, , \, \rangle$ be the standard  
$\mathbb Z$-pairing on $X \times X^\vee$. Then $(X,R,X^\vee,R^\vee)$ is the  
root datum for $\gspin{m}$, with $R$ and $R^\vee$ generated, respectively, by    
\begin{eqnarray*} 
\Delta&=&\{\alpha_1=e_1-e_2,\alpha_2=e_2-e_3,\dots, 
\alpha_{n-1}=e_{n-1}-e_n,\alpha_n=e_n\}, \\  
\Delta^\vee&=&\{\alpha_1^\vee=e_1^*-e_2^*, \alpha_2^\vee=e_2^*-e_3^*,\dots,  
\alpha_{n-1}^\vee=e_{n-1}^*-e_n^*,  
\alpha_n^\vee=2 e_n^*-e_0^*\},  
\end{eqnarray*}  
if $m=2n+1$ and  
\begin{eqnarray*} 
\Delta&=&\{\alpha_1=e_1-e_2,\dots,\alpha_{n-1}=e_{n-1}-e_n,\alpha_n=e_{n-1}+e_n\}, \\ 
\Delta^\vee&=&\{\alpha_1^\vee=e_1^*-e_2^*, \dots, \alpha_{n-1}^\vee=e_{n-1}^*-e_n^*,  
\alpha_n^\vee=e_{n-1}^*+e_n^*-e_0^*\}, 
\end{eqnarray*} if $m=2n$.  
\end{defi} 
 
In the odd case, $\G$ has a Dynkin diagram of type $B_n$: 
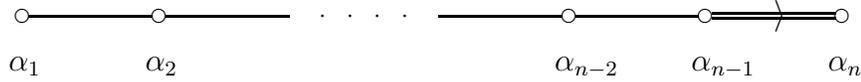
\begin{center} 
\setlength{\unitlength}{2 ex} 
\begin{picture}(35,10) 
\put(0.5,3){\a{1}} 
\put(1,5){\circle{0.5}} 
\put(5.5,3){\a{2}} 
\put(6,5){\circle{0.5}} 
\put(20.5,3){\a{n-2}} 
\put(21,5){\circle{0.5}} 
\put(25.5,3){\a{n-1}} 
\put(26,5){\circle{0.5}} 
\put(30.5,3){\a{n}} 
\put(31,5){\circle{0.5}} 
\thicklines 
\put(1.25,5){\line(1,0){4.5}} 
\put(6.25,5){\line(1,0){4.5}} 
\multiput(12,5)(1,0){4}{\circle*{0.1}} 
\put(16.25,5){\line(1,0){4.5}} 
\put(21.25,5){\line(1,0){4.5}} 
\multiput(26.25,4.9)(0,0.2){2}{\line(1,0){4.5}} 
\put(28.5,4.7){$\rangle$} 
\end{picture} 
\end{center} 
In the even case, it has a Dynkin diagram of type $D_n$: 
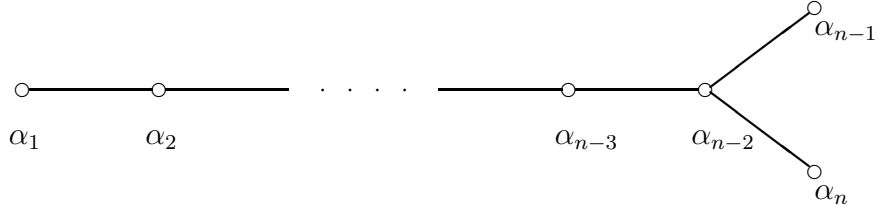
\begin{center} 
\setlength{\unitlength}{2 ex} 
\begin{picture}(35,10) 
\put(0.5,3){\a{1}} 
\put(1,5){\circle{0.5}} 
\put(5.5,3){\a{2}} 
\put(6,5){\circle{0.5}} 
\put(20.5,3){\a{n-3}} 
\put(21,5){\circle{0.5}} 
\put(25.5,3){\a{n-2}} 
\put(26,5){\circle{0.5}} 
\put(30,8){\circle{0.5}} 
\put(30,7){\a{n-1}} 
\put(30,2){\circle{0.5}} 
\put(30,1){\a{n}} 
\thicklines 
\put(1.25,5){\line(1,0){4.5}} 
\put(6.25,5){\line(1,0){4.5}} 
\multiput(12,5)(1,0){4}{\circle*{0.1}} 
\put(16.25,5){\line(1,0){4.5}} 
\put(21.25,5){\line(1,0){4.5}} 
\put(26.2,5.1){\line(4,3){3.6}} 
\put(26.2,4.9){\line(4,-3){3.6}} 
\end{picture} 
\end{center}

\begin{prop}\label{spin-structure} 
The derived group of $\G$ is isomorphic to $\spin{2n+1}$ or $\spin{2n}$,  
the double coverings, as algebraic groups, of special orthogonal groups. In fact,  
$\G$ is isomorphic to  
\[ \frac{\gl{1}\times\spin{m}}{\{(1,1),(-1,c)\}}\, \] 
where $c=\a{n}^\vee(-1)$ if $m=2n+1$ or $c = \a{n-1}^\vee(-1)\a{n}^\vee(-1)$ if  
$m=2n$. The dual of $\G$ is $\gsp{2n}$ if $m=2n+1$ and $\gso{2n}$ if $m=2n$.  
 
Moreover, if $\M$ is the Levi component of a maximal standard parabolic subgroup of $\G$, then  
it is isomorphic to $\gl{k}\times\gspin{m-2k}$ with $k = 1,2,\dots,n$ if $m=2n+1$ and  
$k=1,2,\dots,n-2,n$ if $m=2n$.  
\end{prop} 
\begin{proof} 
See Section 2 of \cite{asgari}. 
\end{proof} 
 
We can also describe the Levi subgroup $\M$ in terms of its root datum.  
Without loss of generality we may assume $\M$ to be maximal. Obviously,  
$\M$ will have the same character and cocharacter lattices as those of  
$\G$. Denote the set of roots of $\M$ by $R_\M$ and its coroots by  
$R_\M^\vee$. They are generated by $\Delta - \{\a{}\}$ and  
$\Delta^\vee - \{\a{}^\vee\}$, respectively, where $\a{} = \a{k}$ unless  
$m=2n$ and $k=n$ in which case $\a{}$ can be either of $\a{n}$ or  
$\a{n-1}$ (resulting in two non-conjugate isomorphic Levi components).  
In the sequel, the case of $m=2n$ and $k=n-1$ is therefore always ruled  
out and we will not repeat this again. 
 
\begin{prop}\label{center-gspin} 
\begin{itemize} 
\item[(a)] The center of $\G$ is given by  
\[ \Z_\G= \left\{ \begin{array} {l@{\quad\mbox{if}\quad}l} 
\A_0 & m=2n+1, \\  
\A_0 \cup (\zeta_0\A_0) & m=2n,\end{array} \right. \]  
where  
\[ \A_0=\left\{e_0^*(\l)\,:\, \l\in\gl{1} \right\} \]  
and $\zeta_0 = e_1^*(-1) e_2^*(-1) \cdots e_n^*(-1)$.  
 
\item[(b)] The center of $\M$ is given by  
\[ \Z_\M= \left\{ \begin{array} {l@{\quad\mbox{if}\quad}l} 
\A_k & m=2n+1, \\   
\A_k \cup (\zeta_k\A_k) & m=2n, \end{array} \right. \]  
where 
\[ \A_k=\left\{e_0^*(\l)e_1^*(\mu)e_2^*(\mu) \cdots e_k^*(\mu)\,:\, \l,  
\mu\in\gl{1} \right\} \]  
and $\zeta_k = e_{k+1}^*(-1) e_{k+2}^*(-1) \cdots e_n^*(-1)$. 
\end{itemize} 
\end{prop} 
\begin{proof} 
The maximal torus $\T$ of $\G$ (or $\M$) consists of elements of the form  
\[ t= \prod_{j=0}^n e_j^*(t_j) \] 
with $t_j\in\gl{1}$. Now $t$ is in the center of $\G$, respectively $\M$,  
if and only if it belongs to the kernel of all simple roots of $\G$,  
respectively $\M$. For $\G$, this leads to  
\[ {t_1}/{t_2} = {t_2}/{t_3} = \cdots = {t_{n-1}}/{t_n}=t_n=1 \]  
if $m=2n+1$ and  
\[ {t_1}/{t_2} = {t_2}/{t_3} = \cdots = {t_{n-1}}/{t_n}=t_{n-1}t_n=1 \]  
if $m=2n$. For $\M$ we get  
\[ {t_1}/{t_2} = {t_2}/{t_3} = \cdots = {t_{k-1}}/{t_k}, \quad  
{t_{k+1}}/{t_{k+2}} = \cdots = {t_{n-1}}/{t_n}=t_n=1 \] if $m=2n+1$ and  
\[ {t_1}/{t_2} = {t_2}/{t_3} = \cdots = {t_{k-1}}/{t_k}, \quad  
{t_{k+1}}/{t_{k+2}} = \cdots = {t_{n-1}}/{t_n}=t_{n-1}t_n=1 \] if $m=2n$.  
These relations prove the proposition. 
\end{proof} 
 
\begin{rem}\label{center} 
When $m=2n$, the non-identity component of $\Z_\G$ can  
also be written as $z^\prime\A_0$, where $z^\prime$ is a non-trivial  
element in the center of $\spin{2n}$, the derived group of $\G$.  
We now specify this element explicitly in terms of the central  
element $z$ of Proposition 2.2 of \cite{asgari}. There is a  
typographical error in the description of $z$ in that article  
which we correct here: 
\[ z= \left\{ \begin{array} {l@{\quad\mbox{if }}l}  
\prod_{j=1}^{n-2}\a{j}^\vee((-1)^j)\cdot\a{n-1}^\vee(-1) 
& n \mbox{ is even,} \\ 
\prod_{j=1}^{n-2}\a{j}^\vee((-1)^j)\cdot 
\a{n-1}^\vee(-\sqrt{-1})\a{n}^\vee(\sqrt{-1}) 
& n \mbox{ is odd.} \end{array} \right. \]  
To compute $z^\prime$ note that with $m=2n$ we have,  
\[ e_1^* + \cdots e_{n-1}^* + e_n^* = \sum_{j=1}^{n-2} j \a{j}^\vee +  
(\frac{n}{2} - 1) \a{n-1}^\vee + \frac{n}{2} \a{n}^\vee + \frac{n}{2} e_0^* \]  
which, when evaluated as a character at $(-1)$, yields  
\[ \zeta_0 = \left\{ \begin{array} {l@{\quad\mbox{if}\quad}l}  
z & n = 4p, \\ 
z e_0^*(\sqrt{-1}) & n = 4p+1, \\ 
c z e_0^*(-1) & n = 4p+2, \\ 
c z e_0^*(-\sqrt{-1}) & n = 4p+3. \end{array} \right. \] 
Therefore, $\zeta_0 \A_0 = z^\prime \A_0$, where $z^\prime$ is an element in  
the center of $\spin{2n}$ given by   
\[ z^\prime = \left\{ \begin{array} {l@{\quad\mbox{if}\quad}l}  
z & n \equiv 0,1,  \mod 4, \\ 
c z & n \equiv 2,3, \mod 4. \end{array} \right. \] 
\end{rem} 
 
\subsection{Root data for $\wgspin{}$ groups}\label{gspin-tilde} 
 
We describe the structure theory for the group $\wgspin{2n}$ as well as  
its standard Levi subgroups in this section. For our future discussion on  
stability of $\g$-factors in Section \ref{stability} we will need to  
work with a group with connected center. However, center of $\gspin{2n}$  
is not connected as we saw in Proposition \ref{center-gspin}. To remedy  
this we define a new group which is just $\gspin{2n}$ with a one-dimensional  
torus attached to it. This group will have a connected center while its  
derived group is the same as that of $\gspin{2n}$, i.e., $\spin{2n}$.  
This will allow us to work with $\wgspin{2n}$ as we will explain  
in Section \ref{stability}. 
 
\begin{defi} Let $\wgspin{2n}$ be the group   
 	\[ \frac{\gl{1}\times\gspin{2n}}{\{(1,1),(-1,\zeta_0)\}}\, \] 
where $\zeta_0$ is as in Proposition \ref{center-gspin}. Note that the  
derived group of $\wgspin{2n}$ is clearly isomorphic to $\spin{2n}$. 
\end{defi} 
 
We now describe the root datum of this group. 
 
\begin{prop}\label{str-wgspin} Let  
\[ X={\mathbb Z}E_{-1}\oplus{\mathbb Z}E_0\oplus{\mathbb Z}E_1 
\oplus\dots\oplus{\mathbb Z}E_n \]  
and  
\[ X^\vee={\mathbb Z}E_{-1}^*\oplus{\mathbb Z}E_0^*\oplus{\mathbb Z}E_1^* 
\oplus\dots\oplus{\mathbb Z}E_n^* \]  
and let $\langle \, , \, \rangle$ be the standard  
$\mathbb Z$-pairing on $X \times X^\vee$. Then $(X,R,X^\vee,R^\vee)$ is the  
root datum for $\wgspin{2n}$, with $R$ and $R^\vee$ generated, respectively, by    
\[\Delta = \{\alpha_1=E_1-E_2,\dots,  
\alpha_{n-1}=E_{n-1}-E_n,\alpha_n=E_{n-1}+E_n-E_{-1}\}, \]  
and  
\[ \Delta^\vee = \{\alpha_1^\vee=E_1^*-E_2^*, \dots,  
\alpha_{n-1}^\vee=E_{n-1}^*-E_n^*, \alpha_n^\vee=E_{n-1}^*+E_n^*-E_0^*\}. \] 
\end{prop} 
 
\begin{proof} 
Our proof will be similar to the proof of Proposition 2.4 of \cite{asgari}.  
We will compute the root datum of $\wgspin{2n}$ using that of $\gspin{2n}$  
described earlier and verify that it can be written as above.  
 
Start with the character lattice of $\gl{1}\times\gspin{2n}$ which can be  
written as the $\zl$-span of $e_0, e_1, \dots, e_n$ and $e_{-1}$. Now,  
characters of $\wgspin{2n}$ are those characters of $\gl{1}\times\gspin{2n}$  
which are trivial when evaluated  
at the element $(-1,\zeta_0)$. Note that $e_i(\zeta_0) = -1$ for $1\le i \le n$,  
$e_0(\zeta_0) =1$, and $e_{-1}(-1) = -1$. This implies that the character  
lattice of $\wgspin{2n}$ can be written as the $\zl$-span of  
$2 e_{-1}, e_0, e_1+e_{-1},\dots,e_n+e_{-1}$. Now, set $E_{-1} = 2 e_{-1}$,  
$E_0 = e_0$ and $E_i = e_{-1} + e_i$ for $1\le i \le n$. Using the  
$\zl$-pairing of the root datum, we can compute a basis for the cocharacter  
lattice which turns out to consist of  
$E_{-1}^* = e_{-1}^*/2 - (e_1^*+\cdots e_n^*)/2$,  
$E_0^* = e_0^*$, and $E_i^* = e_i^*$, for $1\le i \le n$. Writing the  
simple roots and coroots now in terms of the new bases finishes the  
proof. For example,  
\[ \a{n} = e_{n-1} + e_n = (e_{n-1}+ e_{-1}) + (e_n + e_{-1}) - 2 e_{-1}  
= E_{n-1} + E_n - E_{-1} \] and  
\[ \a{n}^\vee = e_{n-1}^* + e_n^* - e_0^* = E_{n-1}^* + E_n^* - E_0^*. \] 
\end{proof} 
 
We can also describe the root datum of any standard Levi subgroup $\M$  
in $\wgspin{2n}$. Again, without loss of generality, we may assume $\M$  
to be maximal. Similar to the case of $\gspin{2n}$, the roots and  
coroots of $\M$ are, respectively, generated by $\Delta - \{\a{k}\}$  
and $\Delta^\vee - \{\a{k}^\vee\}$ for some $k$. The character and  
cocharacter lattices are the same as those of $\wgspin{2n}$. 
 
\begin{prop} 
\begin{itemize} 
\item[(a)] The center of $\wgspin{2n}$ is given by  
\[ \left\{E_0^*(\mu) E_1^*(\l) E_2^*(\l)\cdots E_n^*(\l)  
E_{-1}^*(\l^2)\,:\, \l,\mu\in\gl{1} \right\}, \]  
and is hence connected. 
\item[(b)] The center of $\M$ is given by  
\[ \left\{E_0^*(\mu)E_1^*(\nu)\cdots E_k^*(\nu) E_{k+1}^*(\l)\cdots  
E_n^*(\l) E_{-1}^*(\l^2)\,:\, \l, \mu, \nu \in\gl{1} \right\}, \]  
and is hence connected. 
\end{itemize} 
\end{prop} 
 
\begin{proof} 
The maximal torus of $\wgspin{2n}$ (or $\M$) consists of elements of the form  
\[ t= \prod_{j=-1}^n E_j^*(t_j) \] 
with $t_j\in\gl{1}$. Now $t$ is in the center of $\G$, respectively $\M$,  
if and only if it belongs to the kernel of all simple roots of $\G$,  
respectively $\M$. For $\G$, this leads to  
\[ {t_1}/{t_2} = {t_2}/{t_3} = \cdots = {t_{n-1}}/{t_n}=(t_{n-1}t_n)/t_{-1}=1. \]  
For $\M$ we get  
\[ {t_1}/{t_2} = {t_2}/{t_3} = \cdots = {t_{k-1}}/{t_k}, \quad  
{t_{k+1}}/{t_{k+2}} = \cdots = {t_{n-1}}/{t_n}= (t_{n-1}t_n)/t_{-1}=1. \]  
These relations prove the proposition. 
 
\end{proof} 
 
We also describe the structure of standard Levi subgroups in $\wgspin{2n}$. 
 
\begin{prop}\label{wgspin-levi}  
Standard Levi subgroups of $\wgspin{2n}$ are isomorphic to  
\[  
\gl{k_1} \times \cdots \times \gl{k_r} \times \wgspin{2 l},  
\]  
where $k_1+\cdots+k_r+l=n$.  
\end{prop} 
 
\begin{proof} 
Without loss of generality we may assume $\M$ to be maximal. The character  
and cocharacter lattices of $\M$, which are the same as those of $\G$, were  
described in Proposition \ref{str-wgspin} and can be written as  
\[ (\zl E_1 \oplus \cdots \oplus \zl E_k) \oplus  
(\zl E_{-1} \oplus \zl E_0 \oplus \zl E_{k+1} \cdots \oplus \zl E_n) \]  
and  
\[ (\zl E_1^* \oplus \cdots \oplus \zl E_k^*) \oplus  
(\zl E_{-1}^* \oplus \zl E_0^* \oplus \zl E_{k+1}^* \cdots \oplus \zl E_n^*), \]  
respectively. This along with the description of roots and coroots of  
$\M$ given above implies that the root datum of $\M$ can be written as  
a direct sum of two root data. The first one is now the well-known root  
datum of $\gl{k}$ and the second is just our earlier description of root  
datum of $\wgspin{2(n-k)}$. Therefore, $\M$ is isomorphic to  
$\gl{k} \times \wgspin{2(n-k)}$.  
\end{proof}

\subsection{Root data for $\gsp{2n}$ and $\gso{2n}$}\label{gsp-gso} 
 
We describe the root data for the two groups $\gsp{2n}$ and $\gso{2n}$ in detail.  
Since these two groups are usually introduced as matrix groups, we will also  
describe the root data in terms of their usual matrix representation.  
It will be evident from this description that the two groups $\gspin{2n+1}$  
and $\gsp{2n}$ as well as $\gspin{2n}$ and $\gso{2n}$ are pairs of connected  
reductive algebraic groups with dual root data.  
 
\begin{defi} Let  
\[ X={\mathbb Z}e_0\oplus{\mathbb Z}e_1\oplus\dots\oplus{\mathbb Z}e_n \]  
and  
\[ X^\vee={\mathbb Z}e_0^*\oplus{\mathbb Z}e_1^*\oplus\dots\oplus{\mathbb Z}e_n^* \]  
and let $\langle \, , \, \rangle$ be the standard  
$\mathbb Z$-pairing on $X \times X^\vee$. Then $(X,R,X^\vee,R^\vee)$ is the  
root datum for the connected reductive algebraic group $\gsp{2n}$ or $\gso{2n}$,  
with $R$ and $R^\vee$ generated,  
respectively, by    
\begin{eqnarray*} 
\Delta&=&\{\alpha_1=e_1-e_2,\alpha_2=e_2-e_3,\dots,\alpha_{n-1}=e_{n-1}-e_n, 
\alpha_n= 2 e_n - e_0\}, \\  
\Delta^\vee&=&\{\alpha_1^\vee=e_1^*-e_2^*, \alpha_2^\vee=e_2^*-e_3^*,\dots,  
\alpha_{n-1}^\vee=e_{n-1}^*-e_n^*,  
\alpha_n^\vee=e_n^*\},  
\end{eqnarray*}  
for $\gsp{2n}$ (cf. pages 133--136 of \cite{tadic}) and 
\begin{eqnarray*} 
\Delta&=&\{\alpha_1=e_1-e_2,\dots,\alpha_{n-1}=e_{n-1}-e_n,\alpha_n=e_{n-1}+e_n - e_0\}, \\ 
\Delta^\vee&=&\{\alpha_1^\vee=e_1^*-e_2^*, \dots, \alpha_{n-1}^\vee=e_{n-1}^*-e_n^*,  
\alpha_n^\vee=e_{n-1}^*+e_n^*\}, 
\end{eqnarray*}  
for $\gso{2n}$. 
\end{defi}

The Dynkin diagrams are of type $C_n$ and $B_n$, respectively. A computation similar to  
the proof of Proposition \ref{center-gspin} proves the following:  
 
\begin{prop}\label{cent-gs} 
Let $\G$ be either $\gsp{2n}$ or $\gso{2n}$. Then the center of $\G$  
is given by  
\[ \Z=  
   \left\{e_0^*(\l^2)\,e_1^*(\l)\,\cdots\,e_n^*(\l):\, \l\in\gl{1} \right\} 
\] 
\end{prop} 
 
Alternatively, consider the group defined via  
\[  
\left\{ g\in\gl{2n}\,:\,^tg J g = \mu(g) J \right\}, \] 
where the $2n\times 2n$ matrix $J$ is defined via  
\[ J =  
\left(    \begin{array}{cccccc}  
                                   &   & & & & 1 \\ 
                                      &&&& \adots & \\  
                                   &&& 1&& \\ 
                                   &&-1&&& \\ 
                                   &\adots&&&& \\ 
                                   -1 &&&&&  
                      \end{array}   \right) 
\quad\mbox{or }  
\left(    \begin{array}{cccccc}  
                                   &   & & & & 1 \\ 
                                      &&&& \adots & \\  
                                   &&& 1&& \\ 
                                   &&1&&& \\ 
                                   &\adots&&&& \\ 
                                   1 &&&&&  
                      \end{array}   \right),  
\]  
respectively. The former is the connected reductive algebraic group $\gsp{2n}$.  
However, the latter is not connected as an algebraic group. This group is sometimes  
denoted by $\mbox{GO}_{2n}$ (Section 2 of \cite{ramak}). Its connected component  
is the group $\gso{2n}$ (also denoted by $\mbox{SGO}_{2n}$).   
The maximal split torus in both of these groups can be described as  
\begin{equation}\label{torus-gsp-gso} 
\hT= \left\{ t(a_1,\cdots,a_n,b_n,\cdots,b_1)=\left(    \begin{array}{cccccc}  
                                   a_1&   & & & & \\ 
                                      &\ddots&&&& \\  
                                   &&a_n&&& \\ 
                                   &&&b_n&& \\ 
                                   &&&&\ddots& \\ 
                                   &&&&&b_1  
                      \end{array}   \right) 
                                       : a_i b_i = \mu \right\} 
\end{equation}  
We can now describe $e_i$ and $e_i^*$ in terms of matrices. In either case we have 
 
\begin{eqnarray}\label{e-e*}  
e_0(t) = \mu &, & e_0^*(\l)=t(1,\cdots,1,\l,\cdots,\l) \\ 
e_i(t)= a_i  &,  & e_i^*(\l)=t(1,\cdots,1,\underset{\underset{i}\uparrow}{\l},1, 
\cdots,1,\underset{\underset{2n+1-i}\uparrow}{\l^{-1}},1,\cdots,1),\, 
1\le i \le n. \nonumber 
\end{eqnarray}

\section{Analytic Properties of Local $L$-functions}\label{local} 
 
Let $F$ denote a local field of characteristic zero, either archimedean or  
non-archimedean.  Let $\G_n$ denote the algebraic group $\gspin{2n+1}$ (respectively,  
$\gspin{2n}$) and let $\s$ be an irreducible admissible generic representation  
of $M=\M(F)$ in $G=\G_{r+n}(F)$, where $\M \simeq \gl{r} \times \G_n$ is the Levi subgroup of a  
standard parabolic $\PP$ in $\G_{r+n}$. Let $\hM\simeq\gl{r}(\cx)\times\gsp{2n}(\cx)$  
(respectively, $\hM\simeq\gl{r}(\cx)\times\gso{2n}(\cx)$) be the Levi component  
of the corresponding standard parabolic subgroup $\hP$ in the dual group  
$\hG=\,^LG^0=\gsp{2n}(\cx)$ (respectively, $\hG=\gso{2n}(\cx)$). Let $r$ denote  
the adjoint action of $\hM$ on the Lie algebra of the unipotent radical of $\hP$.  
Then by Proposition 5.6 of \cite{asgari} $r=r_1 \oplus r_2$ if $n \geq 1$  
(respectively, $n \geq 2$) with $r_1=\rho_r\otimes \w{R}$ and $r_2=\mbox{Sym}^2  
\rho_r \otimes \mu^{-1}$ (respectively, $r_2=\wedge^2 \rho_r \otimes \mu^{-1}$).  
Here, $\rho_r$ denotes the standard representation of $\gl{r}(\cx)$, $\w{R}$ denotes  
the contragredient of the standard representation of $\gsp{2n}(\cx)$ (respectively,  
$\gso{2n}(\cx)$), and $\mu$ denotes the multiplicative character defining  
$\gsp{2n}(\cx)$ (respectively, $\gso{2n}(\cx)$). If $n=0$, then $r=r_1$ with  
$r_1=\mbox{Sym}^2 \rho_r \otimes \mu^{-1}$ (respectively, $r_1=\wedge^2 \rho_r  
\otimes \mu^{-1}$). Recall that we have excluded $n=1$ in the even  
case. The Langlands-Shahidi method defines the $L$-functions $L(s,\s,r_i)$  
and $\e$-factors $\e(s,\s,r_i,\psi)$ for $1 \leq i \leq 2$, where $\psi$ is a  
non-trivial additive character of $F$. (In the global setting, it will be the  
local component of our fixed global additive character $\psi$ of Section  
\ref{intro}.) If $\pi$ denotes a representation of $\G_n(F)$ and $\tau$ denotes  
one of $\gl{r}(F)$, then we sometimes employ the following notations for these  
$L$-functions as well as their global analogues: 
 
\begin{eqnarray} 
\label{RS-L} 
L(s, \pi\times\tau) := L(s, \tau\otimes\w{\pi}, \rho_r\otimes\w{R}) =  
L(s,\tau\otimes\w{\pi},r_1), \\ 
\label{RS-ep} 
\e(s, \pi\times\tau,\psi) := \e(s, \tau\otimes\w{\pi}, \rho_r\otimes\w{R},\psi)  
= \e(s,\tau\otimes\w{\pi},r_1,\psi). 
\end{eqnarray} 
 
\begin{prop} \label{conj7.1} Assume that $\s$ is tempered. Then the local $L$-function  
$L(s,\s,r_i)$ is holomorphic for $\Re(s)>0$ for $1 \leq i \leq m$.  
\end{prop} 
\begin{proof} The result is well-known for archimedean $F$. For non-archimedean $F$ this  
is Theorem 5.7 of \cite{asgari}. Here, $i = 1, 2$ and the first $L$-function gives the  
Rankin-Selberg product while the second is twisted symmetric or exterior square. When, $n=0$,  
we only get the second $L$-function. (cf. Proposition 5.6 of \cite{asgari}). 
\end{proof} 
  
\begin{prop}\label{smc} (Standard Module conjecture for $\G_n$) 
Let $\s$ be an irreducible admissible generic representation of $\M(F)$ in $\G_n(F)$  
and let $\nu$ be an element in the positive Weyl chamber. Let $I(\nu, \s)$ be the  
representation unitarily induced from $\nu$ and $\s$, called {\it the standard module},  
and denote by $J(\nu, \s)$ its unique Langlands quotient. Assume that  
$J(\nu, \s)$ is generic. Then, $J(\nu,\s) = I(\nu, \s)$. In particular,  
$I(\nu,\s)$ is irreducible. (A similar result also holds for general linear  
groups \cite{zel}.) 
\end{prop} 
\begin{proof} 
For archimedean $F$, this is due to Vogan for general groups.  
When $F$ is a non-archimedean field, the proof is the subject of W. Kim's thesis  
\cite{kimthesis}, which we rely on.  
 
However, for small values of $n$ we need not rely on \cite{kimthesis} and can  
obtain the result from published results as we now explain: the group $\gspin{5}$  
is isomorphic to $\gsp{4}$ whose derived group is $\sp{4}$. G.~Mui{\'c} has  
proved the Standard Module conjecture for (quasi-split) classical groups  
(Theorem 1.1 of \cite{muic}). The result for $\gspin{5}$ now follows from  
Corollary \ref{sharedDerived} below.  
 
Similarly, note that the derived group of $\gspin{6}$ is isomorphic to  
$\spin{6}\simeq\sl{4}$, hence equal to the derived group of $\gl{4}$. Therefore,  
again by Corollary \ref{sharedDerived}, the result for $\gspin{6}$ follows  
from the Standard Module conjecture for $\gl{6}$.  
\end{proof} 
 
\begin{prop} Let $\G\subset\w{\G}$ be two connected reductive groups whose  
derived groups are equal. Let $\w{\PP}=\w{\M}\N$ be a maximal standard Levi  
subgroup of $\w{\G}$ and let $\PP=\M\N$ be the corresponding one in $\G$ with  
$\M=\w{\M}\cap\G$. Also, let $\w{\T}\subset\w{\M}$ and $\T=\w{\T}\cap\G\subset\M$  
be maximal tori in $\w{\G}$ and $\G$. Let $\w{\s}$ be a quasi-tempered   
representation of $\w{M}=\w{\M}(F)$ and denote by $\s$ its restriction to  
$M=\M(F)$. Write $\s = \underset{i}{\oplus} \s_i$. Let $I(\w{\s})$ denote the  
induced representation $\underset{\w{M}N\uparrow \w{G}}{\mbox{Ind}}  
\w{\s}\otimes 1$ of $\w{G}=\w{\G}(F)$ and $I(\s_i)$ denote  
$\underset{MN\uparrow G}{\mbox{Ind}} \s_i\otimes 1$, a representation of  
$G=\G(F)$. Then, $I(\w{\s})$ is irreducible and standard if and only if  
each $I(\s_i)$ is standard and irreducible. 
\end{prop} 
\begin{proof} 
We only need to address the reducibility questions. Write $\w{\s}|M =  
\oplus_i \s_i$. By irreducibility of $\w{\s}$ and the fact that $\w{M}=\w{T} M$,  
choose \[ \left\{t_1=1,t_2,\dots,t_k\, : \, t_i\in\w{T}=\w{\T}(F)\right\} \]  
such that $\s_i(m)=\s_1(t_i^{-1} m t_i)$. Observe that  
\begin{equation} 
I(\w{\s})|G = \oplus_i I(\s_i). 
\end{equation} 
In fact, if $f_1\in V(\s_1)$, define $f_i(g) = f_1(t_i^{-1} g t_i)$. Then  
$f_i\in V(\s_i)$, the space of $I(\s_i)$, and the representation  
$I(\s_i)(t_i^{-1} g t_i)$ on $V(\s_1)$ is isomorphic to $I(\s_i)$ since  
\begin{equation} 
(I(\s_1)(t_i^{-1} g t_i)f_1)_i = I(\s_i)(g) f_i,  
\end{equation} 
for all $g\in G$. In particular, $I(\s_i)$ is irreducible if and only if  
$I(\s_1)$ is. Observe moreover that $\w{T}$ acts transitively on the set  
of $I(\s_i)$'s using $\w{G}=\w{T}G$.  
 
If each $I(\s_i)$ is irreducible, then $I(\w{\s})$ has to be irreducible.  
In fact, if $(\w{I_1}, \w{V_1})$ is an irreducible subrepresentation of  
$I(\w{\s})$, then  
\begin{equation} 
\w{I_1}|G = \oplus_j I_j, \quad I_j\not=\{0\}, 
\end{equation}  
and given $j$, there exists $i$ such that $I_j\subset I(\s_i)$. Conversely,  
for each $i$ there exists $j$ such that $0 \not= I_j\subset I(\s_i)$. In fact,  
fix $i$ such that $V(\s_i)\cap\w{V_1}\not= 0$. Since $\w{V_1}$ is invariant  
under $\w{T}$, applying $I(\w{\s})(\w{T})$ to this intersection, then implies  
that $V(\s_i)\cap\w{V_j}\not=\{0\}$ for all $i$. Consequently,  
$\{0\}\not= I_j \underset{\not=}\subset I(\s_i)$  which is a contradiction. 
 
Conversely, suppose $I(\w{\s})$ is irreducible but $I(\s_i)$ are (all)  
reducible. Let $V_i$ be an irreducible $G$-subspace of $V(\s_i)$. Then  
\begin{equation} 
\oplus_i I(\w{\s})(t_i) V_i 
\end{equation} 
is a $\w{G}$-invariant subspace of $V(\w{\s})$ which is strictly smaller  
than $V(\w{\s})$, a contradiction. 
\end{proof}  
 
\begin{cor}\label{sharedDerived} 
Suppose $\G$ and $\G^\prime$ are two connected reductive groups having the  
same derived group. Then the standard module conjecture is valid for $G$ if  
and only if it is valid for $G^\prime$.  
\end{cor} 
 
The Langlands-Shahidi method defines the local $L$-functions via the  
theory of intertwining operators. With notation as above, let the  
standard maximal Levi $\M$ in $\G$ correspond to the subset $\theta$  
of the set of simple roots $\Delta$ of $\G$. Then $\theta = \Delta - \{\alpha\}$  
for a simple root $\alpha\in\Delta$. We denote by $w$ the longest  
element in the Weyl group of $\G$ modulo that of $\M$. Then $w$ is the  
unique element with $w(\theta)\subset\Delta$ and $w(\alpha) < 0$. Let  
$A(s, \s, w)$ denote the intertwining operator as in (1.1) on page 278  
of \cite{shahidi:90annals} and let $N(s, \s, w)$ be defined via  
\begin{eqnarray} 
A(s, \s, w) &=& r(s,\s, w) N(s, \s, w), \\    
\label{local-r}  
r(s,\s,w)  &=& \frac{L(s,\s,\w{r}_1) L(2s,\s,\w{r}_2)}{L(1+s,\s,\w{r}_1)  
\e(s,\s,\w{r}_1,\psi) L(1+2s,\s,\w{r}_2) \e(2s,\s,\w{r}_2,\psi)} 
\end{eqnarray}   
 
In fact, the Langlands-Shahidi method inductively defines the $\g$-factors  
using the theory of local intertwining operators out of which the $L$- and  
$\e$-factors are defined via the relation 
\begin{equation} 
\g(s,\s,r_i, \psi) = \e(s,\s,r_i,\psi)\frac{L(1-s,\s,\w{r_i})}{L(s,\s,r_i)}. 
\end{equation} 
 
The following proposition is the main result of this section about analytic  
properties of local $L$-functions which we will use to prove the necessary  
global analytic properties. 
 
\begin{prop}\label{assumptionA} 
Let $\s$ be a local component of a globally generic cuspidal automorphic  
representation of $\M(\Ad)$. Then the normalized local intertwining  
operator $N(s, \s, w)$ is holomorphic and non-zero for $\Re(s)\geq 1/2$.  
\end{prop} 
\begin{proof}  
First assume $\s$ to be tempered. Then by Harish-Chandra  
we know that $A(s, \s, w)$ is holomorphic for $\Re(s)>0$.  
Moreover, for $\Re(s)>0$ we have that $r(s,\s,w)$ is non-zero by definition  
and holomorphic by Proposition \ref{conj7.1}. This implies that $N(s, \s, w)$  
is also holomorphic for $\Re(s)>0$.  
 
Next, assume that $\s$ is not tempered but still unitary. Write $\s = \tau \otimes\w{\pi}$,  
where $\tau$ is a representation of $\gl{r}(F)$ and $\w{\pi}$ is one of $\G_n(F)$.  
(We used $\w{\pi}$ in order to get the usual Rankin-Selberg factors for pairs  
of general linear groups below.) By Proposition \ref{smc}, we can write $\tau$ and  
$\w{\pi}$ as follows: 
\[ \tau = \mbox{Ind} ( \nu^{\a{1}} \tau_1 \otimes \cdots \otimes \nu^{\a{p}} \tau_p \otimes  
\tau_{p+1} \otimes \nu^{-\a{p}} \tau_p \otimes \cdots \otimes \nu^{-\a{1}} \tau_1) \] 
and  
\[ \w{\pi} = \mbox{Ind}( \nu^{\b{1}} \pi_1 \otimes \cdots \otimes {\nu^\b{q}} \pi_q  
\otimes \pi_0), \]  
with $0= \a{p+1}<\a{p}<\cdots<\a{1}<1/2$ and $0<\b{q}<\cdots<\b{1}$ where,  
$\tau_1,\dots,\tau_{p+1}$ and $\pi_1,\dots,\pi_{q}$ are tempered representations of  
the corresponding $\gl{}(F)$, and $\pi_0$ is a generic tempered representation of $\G_t(F)$  
for some $t$. Here, $\nu(\cdot)$ denotes $\big|\det(\cdot) \big|_F$. Therefore,  
\[ N(s,\s, w) = \prod_{\overset{i=1,\dots,p+1}{j=1,\dots,q}} N_1(s\pm\a{i}\pm\b{j})  
\cdot  
\prod_{i=1,\dots,p+1} N_2(s\pm \a{i}), 
\] 
where the terms $N_1(s\pm\a{i}\pm\b{j})$ in the first product are products of  
four rank one operators for $\gl{k} \times \gl{l} \subset \gl{k+l}$ with complex  
parameters $s\pm\a{i}\pm\b{j}$ and the terms $N_2(s\pm \a{i})$ in the second  
product are products of two rank one operators for $\gl{k}\times\G_l\subset\G_{k+l}$  
with complex parameters $s\pm\a{i}$ respectively. 
 
If $\Re(s) \geq 1/2$, then $\Re(s\pm\a{i}) > 0$ for all $i$ and the terms in  
the second product are holomorphic by the first part of  
this proof. The operators in the first product are those associated with  
Rankin-Selberg factors for pairs of general linear groups and, by Lemma 2.10  
of \cite{kim-israel}, they are holomorphic if we show that  
$\Re(s\pm\a{i}\pm\b{j}) > -1$. This happens if $\Re(s-\a{1}-\b{1}) > -1$ or  
$\b{1} < 1$. Therefore, Lemma \ref{b1} below completes the proof of holomorphy part  
of the proposition. 
 
The fact that $N(s,\s,w)$ is a non-vanishing operator now follows from  
applying a result of Y. Zhang (cf. pages 393--394 of \cite{zhang}) to  
our case. Note that in view of Proposition \ref{conj7.1} no assumptions  
are needed in applying \cite{zhang}.  
\end{proof} 
 
\begin{lem}\label{b1} 
If $\pi$ is a supercuspidal representation of $\G_n(F)$ written as above, then  
$\b{1} < 1$.   
\end{lem} 
\begin{proof} 
This is Lemma 3.3 of \cite{kim-israel}. However, note that one should use our $\w{\pi}$  
in the argument. 
\end{proof} 

\section{Stability of $\gamma$-factors}\label{stability} 
 
We continue to denote by $\G_n$ either of the groups $\gspin{2n+1}$ or $\gspin{2n}$  
in this section. In subsections \ref{mel-tras} and \ref{bess-sec} we denote  
by $\wG{n}$ the groups $\gspin{2n+1}$ in the odd case and $\wgspin{2n}$ in  
the even case (see Remark \ref{wgspin-center} below) and $\G$ will denote $\wG{n+1}$  
in either case. 
 
In this section we prove a key local fact, called {\it the stability of $\g$-factors},  
which is what allows us to connect the Langlands-Shahidi $L$- and $\e$-factors  
to those in the converse theorem. Similar results have recently been proved for  
the groups $\so{2n+1}$ in \cite{cogdell-ps,cogdell-kim-ps-sha} and other classical  
groups in \cite{ckpss-classical} which we have followed. A more general result will  
appear in \cite{stability1,stability2}.

Let $F$ denote a non-archimedean local field of characteristic  
zero, i.e, one of $k_v$'s where $v$ is a finite place. Composing a fixed  
splitting with $\psi$ as in \cite{shahidi-imrn}, then defines a generic  
character of $U$ as well as $U_M$ which we still denote by $\psi$.  
Let $\pi$ be an irreducible admissible $\psi$-generic representation of  
$\gspin{2n+1}(F)$ or $\gspin{2n}(F)$, and let $\eta$ be a continuous character  
of $F^\times$. The associated $\g$-factors  
of the Langlands-Shahidi method defined in Theorem 3.5 of \cite{shahidi:90annals}  
will be denoted by $\g(s,\eta\times\pi,\psi)$. They are associated to the  
pair $(\gspin{m+2}, \gl{1}\times\gspin{m})$ of the maximal Levi factor  
$\M = \gl{1}\times\gspin{m}$ in the connected  
reductive group $\gspin{m+2}$, where $m=2n+1$ or $2n$.  
Recall that the $\g$-factor is related to the $L$- and $\e$-factors by  
\begin{equation} 
\g(s,\eta\times,\pi,\psi) = \e(s,\eta\times\pi,\psi)\frac{L(1-s,\eta^{-1}\times 
\w{\pi})}{L(s,\eta\times\pi)}. 
\end{equation} 
 
The main result of this section is the following:  
\begin{thm}\label{stable}  
Let $\pi_1$ and $\pi_2$ be irreducible admissible generic  
representations of $\gspin{m}(F)$ with equal central characters  
$\oo_{\pi_1}=\oo_{\pi_2}$. Then for every  
sufficiently ramified character $\eta$ of $F^\times$ we have  
\[ \g(s,\eta\times\pi_1,\psi) = \g(s,\eta\times\pi_2,\psi). \] 
\end{thm} 
 
The proof of this theorem is the subject matter of this section including a  
review of some facts about {\it partial Bessel functions}.  
 
\subsection{$\g(s,\eta\times\pi)$ as Mellin transform of a Bessel function}\label{mel-tras} 
 
Recall that $\wG{n}$ denotes either $\gspin{2n+1}$ or $\wgspin{2n}$ and $\G=\wG{n+1}$.  
We will refer to $\G=\gspin{2n+3}$ as {\it odd} and $\G=\wgspin{2n+2}$ as {\it even}  
in the rest of this section.  
 
\begin{rem}\label{wgspin-center}  
We will need to assume that the center of our group $\G$ is  
connected for proof of Proposition \ref{asymp1} below. This is not  
true if $\G$ is taken to be $\gspin{2n+2}$ as we pointed out in Proposition \ref{center-gspin}.  
To remedy this we can alternatively work with the group $\wgspin{2n+2}$  
of Section \ref{gspin-tilde}. Since $\wgspin{2n+2}$ has the same  
derived group as $\gspin{2n+2}$ its corresponding $\g$-factors  
are the same as those of $\gspin{2n+2}$ since they (and, in fact,  
the local coefficients via which they are defined) only depend  
on the derived group of our group. This has no effect on the arguments  
of the next few subsections as all of our crucial  
computations take place inside the derived group. 
\end{rem}

Let $\G$ be as above with a fixed Borel subgroup $\B=\T\U$ as before. We  
continue to denote its root data by $(X,R,X^\vee,R^\vee)$ which we have  
described in detail in Section \ref{str}.  
Consider the maximal parabolic subgroup $\PP=\M\N$ in $\G$, where $\N\subset\U$  
and the Levi component, $\M$, is isomorphic to $\gl{1}\times\G_n$. The standard  
Levi subgroup $\M\supset\T$ corresponds to the subset $\theta = \Delta - \{\a{1}\}$  
of the set of simple roots $\Delta=\{\a{1},\a{2},\dots,\a{n},\a{n+1}\}$ of  
$(\G,\T)$ with notation as in Section \ref{str}. Let $w$ denote the unique  
element of the Weyl group of $\G$ such that $w(\theta)\subset\Delta$ and  
$w(\a{1})<0$. Notice that the parabolic subgroup $\PP$ is self-associate, i.e.,  
$w(\theta)=\theta$. We will denote by $G$, $P$, $M$, $N$, $B$, etc., the  
groups $\G(F)$, $\PP(F)$, $\M(F)$, $\N(F)$, $\B(F)$, etc., in what follows.  
Also denote the opposite parabolic subgroup to $P$ by $\bP=M\bN$ .  
 
Let $\Z=\Z_{\G}$ and $\Z_{\M}$ be the centers of $\G$ and $\M$ respectively. The following  
is assumption 5.1 of \cite{shahidi-imrn} for our cases. We will need this when dealing  
with Bessel functions. 
 
\begin{prop}\label{Assumption5.1} 
There exits an injection $e^*\,:\,F^\times\longrightarrow  
Z_G\backslash Z_M$ such that for all $t\in F^\times$ we have  
$\a{1}(e^*(t))=t$. 
\end{prop} 
\begin{proof} 
We define $e^*(t)$ to be the image in $Z_G\backslash Z_M$ of  
$e_1^*(t)$ in the odd case and that of $E_1^*(t)$ in the even case. 
The proposition is now clear from our explicit descriptions in  
Section \ref{str}.  
\end{proof} 
 
Denote the image of $e^*$ by $Z_M^0$ as in \cite{shahidi-imrn}. (Note  
that \cite{shahidi-imrn} uses the notation $\aa^\vee$ for $e^*$.)  
 
We now review some standard facts about the reductive group $G$  
whose proofs could be found in either of \cite{springer} or  
\cite{steinberg}, for example. For $\aa\in R$ let  
$u_{\aa}\,:\,F\longrightarrow G$ be the root group homomorphism  
determined by the equation  
\begin{equation}\label{RootGroup} 
t u_{\aa}(x) t^{-1} = u_{\aa}(\aa(t) x) \quad t\in T, x\in F.  
\end{equation}  
Then $u_{\aa}(x)$ is additive in the variable $x$.   
Moreover, define $w_{\aa}\,:\,F^\times\longrightarrow G$ by  
\begin{equation}\label{general_w} 
w_{\aa}(\lambda) = u_{\aa}(\lambda)u_{-\aa}(-\lambda^{-1})u_{\aa}(\lambda). 
\end{equation}  
Also set $w_{\aa} := w_{\aa}(1)$. 
Then,  
\begin{equation}\label{w_alpha} 
w_{\aa}(\lambda)=\aa^\vee(\lambda) w_{\aa} =w_{\aa}\, \aa^\vee(\lambda^{-1}),  
\end{equation} 
where $\aa^\vee$ is the coroot corresponding to $\aa$.  
The element $w_\aa$ normalizes $T$ and we will denote 
its image in the Weyl group by $\w{w}_\aa$.
 
\begin{rem}\label{u_alpha} 
Our choice of $w_{\aa}$ is indeed the same as $n_{\aa}$ on  
page 133 of \cite{springer}. 
This choice differs up to a sign from those made in (4.43),  
(4.19) or (4.56) of \cite{shahidi-imrn} requiring $w_\aa$  
to be the image of  
$\left(\begin{matrix}0&-1\\1&0\end{matrix}\right)$ under  
the homomorphism from $\sl{2}$ into $\G$ determined by $\aa$.   
The latter choice would introduce an occasional negative sign in some of  
the equations, e.g. (\ref{w_alpha}) above. Of course, this   
choice is irrelevant to final results and we have chosen Springer's  
since we will be using some detailed information on structure  
constants from \cite{springer} in what follows (e.g., Lemma \ref{w_gamma}).  
\end{rem} 
 
Recall that  
\begin{eqnarray}  
w_{-\aa}(\lambda)&=&w_{\aa}(-1/\lambda), \\  
w_\aa ^2 &=& \aa^\vee(-1), \label{w_alpha^2}\\ 
w_{-\aa} &=& w_\aa^{-1}. 
\end{eqnarray}  
 
For any two linearly independent roots $\aa$ and $\bb$ in $R$  
and a total order on $R$, which we fix now, we have  
\begin{equation}\label{comm}  
u_{\aa}(x)u_{\bb}(y)u_{\aa}(-x) = u_{\bb}(y)   
\prod_{\overset{i,j>0}{i\aa+j\bb\in R}} u_{i\aa+j\bb}(c_{ij} x^i y^j)  
\end{equation} 
for certain {\it structure constants} $c_{ij} = c_{\aa,\bb;i,j}$.  
In particular, if there are no roots of the form $i\aa+j\bb$  
with $i,j>0$, then  
$u_{\aa}(x)u_{\bb}(y)=u_{\bb}(y)u_{\aa}(x)$.

We recall the following result. 
 
\begin{prop}\label{w_s_alpha}  
Let $\aa$ and $\bb$ be two arbitrary linearly independent roots and let  
$(\bb-c\aa,\cdots,\bb+b\aa)$ be the $\aa$-string through $\bb$.  
Then,  
\begin{equation}\label{conjugate_w_beta} 
w_{\aa} w_{\bb}(x) w_{\aa}^{-1}  =  w_{\w{w}_\aa(\bb)} (d_{\aa,\bb} x),  
\end{equation} 
where  
\begin{equation}\label{d_alpha}  
 d_{\aa,\bb} =  
\sum_{i=\max(0,c-b)}^c (-1)^i \,  
c_{-\aa,\bb;i,1} \, c_{\aa,\bb-i\aa;i+b-c,1}. 
\end{equation} 
Moreover, 
\begin{equation}\label{minus_a}
d_{-\aa,\bb} = (-1)^{\langle \bb, \aa^\vee\rangle} d_{\aa,\bb} 
\end{equation}
and 
\begin{equation}\label{d_w_a}
d_{\aa,\bb} d_{\aa,\w{w}_\aa(\bb)} = (-1)^{\langle \bb, \aa^\vee\rangle}. 
\end{equation}

\end{prop} 
\begin{proof} 
This is Lemma 9.2.2 of \cite{springer}. Note that Springer  
defines $d_{\aa,\bb}$ via  
\[ w_\aa u_\bb(x) w_\aa^{-1} = u_{\w{w}_\aa(\bb)}(d_{\aa,\bb} x) \]  
which, using (\ref{general_w}), immediately implies (\ref{conjugate_w_beta}).   
\end{proof} 
 
Denoting the image of $u_{\aa}$ in $G$ by $U_{\aa}$ notice that  
$M$ is generated by $T$ and $U_{\aa}$'s with $\aa$ ranging over $\Sigma(\theta)$,  
the set of all (positive and negative) roots spanned by  
$\a{2},\a{3},\dots,\a{n},\a{n+1}$, while $N$ is generated by  
$U_{\aa}$'s where $\aa$ ranges over $R(N)=R^+ - \Sigma(\theta)$,  
the set of positive roots of $G$ not in $M$ (i.e., involving a positive  
coefficient of $\a{1}$ when written as a sum of simple roots with  
non-negative coefficients) and $\bN$ is generated by $U_{\aa}$'s  
where $\aa$ ranges over $R(\bN)=R^- - \Sigma(\theta)$,  
the set of negative root of $G$ not in $M$  
(i.e., involving a negative coefficient of $\a{1}$ when written  
as a sum of simple roots with non-positive coefficients). Let  
$U_M=U \cap M$. Then $U_M$ is generated  
by $U_{\aa}$'s with $\aa \in \Sigma(\theta)^+ = \Sigma(\theta) \cap R^+$.  
 
The group $M$ acts via the adjoint action on $N$; in particular,  
both $U_M$ and $Z_M^0$ act on $N$. We are interested in the orbits  
of the adjoint action of $Z_M^0 U_M$ on $N$.  
 
\begin{lem}\label{U_M} 
Up to a subset of measure zero of $N$ the following is a complete set of  
representatives for the orbits of $N$ under the adjoint action of $U_M$. 
\[ U_M \backslash N \simeq  
\left\{u_{\a{1}}(a) u_{\g}(x)\,:\, a\in F^\times,\, x\in F  \right\}, \]  
where  
\[ \g =  
\begin{cases} 
\a{1} + 2 \a{2} + \cdots + 2 \a{n+1}  
& \mbox{ if }\, \G \mbox{ is odd, } \\ 
\a{1} + 2 \a{2} + \cdots + 2 \a{n-1} + \a{n} + \a{n+1} 
& \mbox{ if }\, \G \mbox{ is even, }  
\end{cases} 
\] 
is the longest positive root in $G$. 
\end{lem}  
\begin{proof} 
Using the same Bourbaki notation as in Section \ref{str},  
$R(N)$ is given by (\cite{bourbaki}) 
 
\begin{equation}\label{beta-odd}   
\left\{  \begin{array}{l} 
	\a{1}, \a{1}+\a{2}, \dots, \a{1}+\a{2}\cdots+\a{n+1},  
	\a{1}+\a{2}+\cdots+\a{n}+2\a{n+1}, \\  
 
	\a{1}+\a{2}+\cdots+\a{n-1}+2\a{n}+2\a{n+1}, \cdots,  
	\g=\a{1}+2\a{2}+\cdots+\\ 2\a{n}+2\a{n+1}  
         \end{array} \right\}  
\end{equation} 
for the odd case and  
\begin{equation}\label{beta-even}    
\left\{  \begin{array}{l}  
	\a{1}, \a{1}+\a{2}, \dots, \a{1}+\a{2}\cdots+\a{n-1}+\a{n},  
	\a{1}+\a{2}\cdots+\a{n-1}+\\ 
 
	\a{n+1}, \a{1}+\a{2}\cdots+\a{n-1}+\a{n}+\a{n+1}, 
	\a{1}+\a{2}+\cdots+\a{n-2}+\\ 2\a{n-1}+\a{n}+\a{n+1},  
	\a{1}+\a{2}+\cdots+2\a{n-2}+2\a{n-1}+\a{n}+\a{n+1}, \\  
 
	\cdots, \g=\a{1}+2\a{2}+\cdots+2\a{n-1}+\a{n}+\a{n+1}  
         \end{array} \right\}  
\end{equation} for the even case.

An arbitrary element $n\in N$ is of the form  
\begin{equation}\label{n-expr} 
n = \prod_{\aa\in R(N)} u_{\aa} (x_\aa)  
\end{equation} 
with $x_\aa \in F$. The ordering in the product can be arbitrarily chosen  
since any linear combination with positive integer coefficients of two  
roots in $R(N)$ would have $\a{1}$ with an integer coefficient of at  
least two which can not be a root, hence by (\ref{comm}) any two terms  
in the above product commute. We make use of this fact in the rest of  
the proof. 
 
Observe that the set $R(N)$ has the property that if $\aa$ belongs to  
$R(N)$, then so does $\g-\aa+\a{1}$. Notice that if   
$\aa^\prime\in R(N) - \left\{\a{1},\g \right\}$,  
then $\bb = \aa^\prime - \a{1}\in \Sigma(\theta)$  
and $\bb > 0$, hence $g=u_\bb(x) \in U_M$ for any $x \in F$. Fix one  
such $\bb$ and consider the adjoint action of $g$ on $n$: 
\[ g n g^{-1} = \prod_{\aa\in R(N)} g u_{\aa}(x_\aa) g^{-1} =  
\prod_{\aa\in R(N)} u_\bb(x) u_{\aa}(x_\aa) u_\bb(-x).  \] 
We now look at each term in this product: if $i\bb+j\aa\not\in R$  
for positive $i$ and $j$, then by (\ref{comm}) the term is equal to  
$u_{\aa}(x_\aa)$. This is the case for all $\aa\in R(N)$ unless  
$\aa=\a{1}$ in which case $\bb+\aa=\aa^\prime$ is a root or  
$\aa=\g-\bb = \g-\aa^\prime+\a{1}$ (which does belong to $R(N)$ by  
the above observation) in which case $\bb+\aa=\g$ is again a root.  
Therefore,  
\[ \prod_{\aa\in R(N)} u_\bb(x) u_{\aa}(x_\aa) u_\bb(-x) =  
\prod_{\aa\in R(N)} u_{\aa}(y_\aa), \]  
where  
\[ y_\aa = \begin{cases} 
x_\aa + C x_{\a{1}} x_\bb\quad&\mbox{ if } \aa = \aa^\prime \\ 
x_\aa + C^\prime x_{\g-\aa^\prime+\a{1}}x_\bb\quad&\mbox{ if } \aa=\g \\ 
x_\aa\quad&\mbox{ otherwise. } 
\end{cases} \] 
Here, $C, C^\prime\in F^\times$ are the appropriate structure constants  
as in (\ref{comm}). 
 
Assuming $x_\a{1} \not= 0$, which only excludes a subset of $n \in N$ of  
measure zero, we can choose $x_\bb\in F$ appropriately in order to  
have $x_\aa + C x_\a{1} x_\bb = 0$. Applying this process for all the  
$\bb$ in $\Sigma(\theta)$ described above we can make all $x_\aa$  
in (\ref{n-expr}) equal to zero except for $x_\a{1}$ and $x_\g$. In  
the process the value of $x_\a{1}$ does not change but the value of  
$x_\g$ may change. We let $a=x_\a{1}$ and let $x$ be the final value of  
$x_\g$. This proves the lemma. 
\end{proof} 
 
We now consider the adjoint action of $Z_M^0$.  
 
\begin{lem}\label{Z_M^0} 
Let $n=u_{\a{1}}(a) u_\g(x) \in N$ with $a \in F^\times$ and $x \in F$. Then,  
$n$ and $u_{\a{1}}(1) u_\g(y)$ are in the same conjugacy class of the adjoint  
action of $Z_M^0$ on $N$ for some $y \in F$. 
\end{lem} 
\begin{proof} 
For $z = e_1^*(\lambda) \in Z_M^0$, in the odd case, we have  
\begin{eqnarray*} 
z n z^{-1} &=& e_1^*(\lambda) u_{\a{1}}(a) e_1^*(\lambda^{-1}) \,\,   
e_1^*(\lambda) u_\g(x) e_1^*(\lambda^{-1}) \\ 
&=& u_{\a{1}}(\a{1}(e_1^*(\lambda)) \, a) \,\, u_\g(\a{1}(e_1^*(\lambda)) \, x) \\ 
&=& u_{\a{1}}(\lambda^{\langle \a{1}, e_1^*\rangle} \, a) \,\, 
  u_\g(\lambda^{\langle \g, e_1^*\rangle} \, x) \\ 
&=& u_{\a{1}}(\lambda \, a) \,\, u_\g(\lambda \, x). 
\end{eqnarray*} 
In the even case $e_1^*$ above should be replaced by $E_1^*$.  
Take $\lambda=1/a$ and $y=x/a$ to finish the proof. 
\end{proof} 
 
Lemmas \ref{U_M} and \ref{Z_M^0} immediately imply the following:  
 
\begin{cor}\label{conjugacy} 
Up to a subset of measure zero of $N$ the following is a complete set of  
representatives for the orbits of $N$ under the adjoint action of $Z_M^0 U_M$:  
\[ Z_M^0 U_M \backslash N \simeq  
\left\{u_{\a{1}}(1) u_{\g}(x)\,:\, x\in F  \right\}. \]  
\end{cor} 
 
If we set  
\begin{equation}\label{w0} 
w_0 =  \begin{cases} 
w_{\a{1}}w_{\a{2}}\cdots w_{\a{n}}w_{\a{n+1}}w_{\a{n}}\cdots w_{\a{2}}w_{\a{1}} 
& \mbox{ if }\, \G \mbox{ is odd, } \\ 
w_{\a{1}}w_{\a{2}}\cdots w_{\a{n-1}} 
w_{\aa\aa}w_{\a{n-1}}\cdots w_{\a{2}}w_{\a{1}}  
& \mbox{ if }\, \G \mbox{ is even. }  
\end{cases} 
\end{equation} 
where $w_{\aa\aa}$ is the product of the commuting elements  
$w_{\a{n}}w_{\a{n+1}}=w_{\a{n+1}}w_{\a{n}}$, then $w_0$ is the  
representative in $G$ of the Weyl group element $w$  
mentioned earlier.  
 
\begin{rem}\label{waa} 
Note that in \cite{shahidi-imrn}  
the analogue of our element $w_{\aa\aa}$ is denoted by  
$w_{\a{n+1}}$ as explained in (4.45) of that article. However, our  
current choices, which applies equally well to $\so{2n}$ or other cases  
treated in Section 4 of \cite{shahidi-imrn}, would  
replace the two commuting matrices on the left hand side of (4.45)  
of \cite{shahidi-imrn} with their transposes (cf. Remark \ref{u_alpha}). 
\end{rem}

As in \cite{shahidi-imrn} we are interested in elements $n \in N$ such that  
\begin{equation}\label{decomp} 
w_0^{-1} n = m n^\prime \bn \in P\bN. 
\end{equation}

The decomposition in (\ref{decomp}) is clearly unique and we would  
like to compute the  
$m-$, $n^\prime-$, and $\bn-$parts of an element $n$ as in  
Corollary \ref{conjugacy}. We will do this in Proposition \ref{conj-reps}.  
First we prove the following auxiliary lemma.   
 
\begin{lem}\label{w_gamma} 
We can normalize the $u_\aa$'s such that the element $w_\gamma$ satisfies     
\[ \g^\vee(d) w_\g = w_\g (d) =   
\begin{cases} 
\begin{matrix} 
w_{\a{2}}\cdots w_{\a{n}}w_{\a{n+1}}w_{\a{n}}\cdots w_{\a{2}} w_{\a{1}} \\ 
\cdot \, 
w_{\a{2}}^{-1}\cdots w_{\a{n}}^{-1}w_{\a{n+1}}^{-1}w_{\a{n}}^{-1}\cdots w_{\a{2}}^{-1} 
\end{matrix} 
& \mbox{ if }\, \G \mbox{ is odd, } \\ 
 & \\  
\begin{matrix} 
w_{\a{2}}\cdots w_{\a{n-1}} w_{\aa\aa} w_{\a{n-1}}\cdots w_{\a{2}} w_{\a{1}} \\ 
\cdot \, w_{\a{2}}^{-1}\cdots w_{\a{n-1}}^{-1} w_{\aa\aa}^{-1} w_{\a{n-1}}^{-1}  
\cdots w_{\a{2}}^{-1}  
\end{matrix} 
& \mbox{ if }\, \G \mbox{ is even, }  
\end{cases} 
\] 
where 
\begin{equation}\label{ddd}
d = \begin{cases}
(-1)^n & \mbox{ if }\, \G \mbox{ is odd, } \\
(-1)^{n-1} & \mbox{ if }\, \G \mbox{ is even. }
\end{cases}
\end{equation}
\end{lem} 

\begin{rem}
The above $d$ in the odd case is slightly different from the corresponding value 
for the group $\so{2n+3}$ carried out in 4.2.1 of \cite{ckpss-classical}, i.e., 
it differs by a factor of $-1$. The reason for this discrepancy is that the 
representative we have fixed above for the longest element of the Weyl group 
would, in the case of the group $\so{3}$, lead to 
$\left(\begin{matrix}
 & & -1 \\ & -1 & \\ -1 & &  
\end{matrix}\right).$ 
This is the correct element which should have been used in 
4.2.1 of \cite{ckpss-classical} and \cite{shahidi-imrn} instead of 
$\left(\begin{matrix}
 & & 1 \\ & -1 & \\ 1 & &  
\end{matrix}\right)$. 
The latter is the corresponding Weyl group representative 
for the group $\sl{3}$. However, since the group $\so{3}$ does not 
have a Cartan of the same dimension as that of $\sl{3}$ there is no natural 
way (i.e., not requiring a choice of basis) of embedding it in $\sl{3}$. Therefore, 
there is no reason why the representative for $\so{3}$ would be the same as 
that of $\sl{3}$. Of course, both of these two matrices correspond to the same 
Weyl group element since they only differ by a diagonal matrix in $\so{3}$. 
In the notation of the present paper we can fix $u_\aa(\ )$ and $u_{-\aa}(\ )$ in $\so{3}$ 
such that 
$w_\aa({\lambda}) =\left(\begin{matrix}
 & & -\lambda^2 \\ 
 & -1 & \\ 
 -1/\lambda^2 & &  \end{matrix}\right)$. 
The choice of representative above is now simply $w_\aa(1)$. 
\end{rem}
 
\begin{proof} 
We begin by noting that  
\begin{equation}\label{RootGamma} 
\g =  
\begin{cases} 
\w{w}_{\a{2}}\cdots\w{w}_{\a{n}}\w{w}_{\a{n+1}}\w{w}_{\a{n}}\cdots\w{w}_{\a{2}}(\a{1})  
& \mbox{ if }\, \G \mbox{ is odd, } \\ 
\w{w}_{\a{2}}\cdots\w{w}_{\a{n-1}}\w{w}_{\a{n}}\w{w}_{\a{n+1}}\w{w}_{\a{n-1}}  
\cdots \w{w}_{\a{2}}(\a{1})   
& \mbox{ if }\, \G \mbox{ is even. }  
\end{cases} 
\end{equation} 
Let $\b{1}=\a{1}$ and denote by $\b{i}$ the consecutive images of $\b{1}$ under the  
first $i-1$ Weyl group elements above for  
\[ \begin{cases}  
1 \le i \leq 2n  
& \mbox{ if }\, \G \mbox{ is odd, } \\  
1 \le i \leq 2n-1 
& \mbox{ if }\, \G \mbox{ is even. }  
\end{cases} 
\] 
In fact, the $\b{i}$'s are precisely the roots listed in (\ref{beta-odd}) 
and (\ref{beta-even}) and in the same order. 

Now apply (\ref{conjugate_w_beta}) repeatedly to conclude that 
the right hand side of the expression of the statement of the Lemma is 
equal to $ w_\g(d)$, where 
\begin{equation*} 
d =  
\begin{cases}  
d_{\a{2},\b{1}} \cdots d_{\a{n},\b{n-1}} d_{\a{n+1},\b{n}} d_{\a{n},\b{n+1}} 
\cdots d_{\a{2},\b{2n-1}}
& \mbox{ if }\, \G \mbox{ is odd, } \\ 
 & \\  
d_{\a{2},\b{1}} \cdots d_{\a{n-1},\b{n-2}} d_{\a{n},\b{n-1}} d_{\a{n+1},\b{n}} 
\cdots  d_{\a{2},\b{2n-2}} 
& \mbox{ if }\, \G \mbox{ is even. } \\ 
\end{cases} 
\end{equation*}  

For $1 \leq i \leq n-1$ in the odd case and for $1 \leq i \leq n$ in the  
even case, the $\a{i+1}$-string through $\b{i}$ is  
$(\b{i},\b{i}+\a{i+1})$, i.e., $c=0$ and $b=1$ in the notation of  
Proposition \ref{w_s_alpha}. In the odd case with $i=n$, the  
$\a{i+1}$-string through $\b{i}$ is $(\b{i},\b{i}+\a{i+1},\b{i}+2\a{i+1})$,  
i.e., $c=0$ and $b=2$. Similarly, for $1 \leq j \leq n-1$ in the 
odd case the $\a{n-j-1}$-string  
through $\b{n+j}$ is $(\b{n+j},\b{n+j}+\a{n-j-1})$, i.e., $c=0$ and $b=1$.  
Also,  for $1 \leq j \leq n-2$ in the even case the $\a{n-j}$-string  
through $\b{n+j}$ is $(\b{n+j},\b{n+j}+\a{n-j})$, i.e., $c=0$ and $b=1$. 
Putting all these together and using (\ref{d_alpha}) we can write 
\[ \begin{array}{lcl}
d_{\a{2},\b{1}} &=& c_{\a{2},\b{1};1,1} \\
&\vdots& \\
d_{\a{n},\b{n-1}} &=& c_{\a{n},\b{n-1};1,1} \\
d_{\a{n+1},\b{n}} &=& c_{\a{n+1},\b{n};2,1} \\
d_{\a{n},\b{n+1}} &=& c_{\a{n},\b{n+1};1,1} \\
&\vdots& \\
d_{\a{2},\b{2n-1}} &=& c_{\a{2},\b{2n-1};1,1} \\
\end{array}
\quad\mbox{ and }\quad
\begin{array}{lcl}
d_{\a{2},\b{1}} &=&c_{\a{2},\b{1};1,1} \\  
&\vdots& \\
d_{\a{n-1},\b{n-2}} &=& c_{\a{n-1},\b{n-2};1,1} \\ 
d_{\a{n},\b{n-1}}  &=&  c_{\a{n},\b{n-1};1,1} \\
d_{\a{n+1},\b{n}} &=&  c_{\a{n+1},\b{n};1,1} \\ 
&\vdots& \\ 
d_{\a{2},\b{2n-2}} &=&  c_{\a{2},\b{2n-2};1,1} \\
\end{array}
\]
in the odd and even cases respectively.

We can now normalize the $u_\aa$'s such that we have 
$c_{\a{i},\a{i+1};1,1} = 1$ and, in the odd case, 
$c_{\a{n},\a{n+1};1,2} = -1$. These normalizations 
are motivated by the explicity matrix realizations of 
the related groups $\so{2n+3}$ and $\so{2n+2}$ such 
as those in \cite{shahidi-imrn}. The values of other 
structure constants are now uniquely determined by 
these (cf. Lemma 9.2.3 of \cite{springer}). We then 
get 
\[ \begin{array}{lcl}
c_{\a{2},\b{1};1,1} &=& -1 \\
&\vdots& \\
c_{\a{n},\b{n-1};1,1} &=& -1\\
c_{\a{n+1},\b{n};2,1} &=& -1\\
c_{\a{n},\b{n+1};1,1} &=& +1\\
&\vdots& \\
c_{\a{2},\b{2n-1};1,1} &=& +1\\
\end{array}
\quad\mbox{ and }\quad
\begin{array}{lcl}
c_{\a{2},\b{1};1,1} &=& -1\\  
&\vdots& \\
c_{\a{n-1},\b{n-2};1,1} &=& -1\\ 
c_{\a{n},\b{n-1};1,1} &=& -1\\
c_{\a{n+1},\b{n};1,1} &=& +1\\ 
&\vdots& \\ 
c_{\a{2},\b{2n-2};1,1} &=& +1\\
\end{array}
\]
Therefore, in the expression for $d$ the first $n$ terms in 
the odd case and the first $n-1$ terms in the even case are 
equal to $-1$ and others are equal to $+1$. This implies that 
$d=(-1)^{n}$ in the odd case and $d=(-1)^{n-1}$ in the even case.
\end{proof}

\begin{prop}\label{conj-reps} 
Assume that $x \in F$  and $n=u_\a{1}(1) u_\g(x)$ satisfies (\ref{decomp}).  
Moreover, assume that $x$ is non-zero (which rules out only a subset  
of $N$ of measure zero). 
Then,  
\begin{eqnarray*} 
m &=& w^\prime \g^\vee(d/x)		\\  
n^\prime &=& u_\g(-x) u_\a{1}(-1)	\\ 
\bn &=& u_{-\g}(1/x) u_{-\a{1}}(1),	 
\end{eqnarray*} 
where $d$ is as in (\ref{ddd}) and  
\[  
w^\prime =  
	\begin{cases} 
w_{\a{2}}^{-1} \cdots w_{\a{n}}^{-1}w_{\a{n+1}}^{-1}w_{\a{n}}^{-1}\cdots w_{\a{2}}^{-1}  
& \mbox{ if }\, \G \mbox{ is odd, } \\ 
w_{\a{2}}^{-1}\cdots w_{\a{n-1}}^{-1}w_{\aa\aa}^{-1}w_{\a{n-1}}^{-1}\cdots w_{\a{2}}^{-1} 
& \mbox{ if }\, \G \mbox{ is even, }  
	\end{cases} 
\]			 
with $w_{\aa\aa}$ again as in Remark \ref{waa}. 
Moreover, we could also write $m$ as  
\begin{equation}\label{m_alternate} 
m = \a{1}^\vee(d/x) w^\prime,  
\end{equation} 
which is analogous to Propositions 4.4 and 4.8 of \cite{shahidi-imrn} 
modulo our Remark \ref{u_alpha}.  
\end{prop} 
 
\begin{proof} 
By uniqueness of decomposition in (\ref{decomp}) it is enough to prove that the  
above values satisfy (\ref{decomp}). This is a straight forward  
computation utilizing (\ref{comm}) multiple times. 
 
First, observe that if $i,j>0$ are integers, then $i \a{1} + j \g$  
can not be a root. Hence, $u_\g(\cdot)$ and $u_{\a{1}}(\cdot)$ commute  
by (\ref{comm}). Also, $i \a{1} + j (-\g)$ can not  
be a root, which again implies by (\ref{comm}) that  
$u_{-\g}(\cdot)$ and $u_{\a{1}}(\cdot)$ commute. Similarly,  
$u_{\g}(\cdot)$ and $u_{-\a{1}}(\cdot)$ commute. 

Moreover, by (\ref{w0}) we have 
\begin{equation} 
w_0 = w_{\a{1}} {w^\prime}^{-1} w_{\a{1}}
\end{equation} 
and by Lemma \ref{w_gamma}, 
\begin{equation} 
w_\g(d) = {w^\prime}^{-1} w_{\a{1}} w^\prime. 
\end{equation}
Now, 
\begin{eqnarray*}
w_0^{-1} n &=& w_{\a{1}}^{-1} {w^\prime} w_{\a{1}}^{-1} 
u_{\a{1}}(1) u_\g(x)\\
&=& w_{\a{1}}^{-1} {w^\prime} u_{\a{1}}(-1) u_{-\a{1}}(1) 
u_{\a{1}}(-1) u_{\a{1}}(1) u_\g(x)\\ 
&=& w_{\a{1}}^{-1} {w^\prime} u_\g(x) u_{\a{1}}(-1)  u_{-\a{1}}(1) \\ 
&=& w_{\a{1}}^{-1} {w^\prime} w_\g(x) u_\g(-x) 
u_{-\g}(1/x) u_{\a{1}}(-1)  u_{-\a{1}}(1)\\
&=& w_{\a{1}}^{-1} {w^\prime} w_\g(x) \cdot u_\g(-x) u_{\a{1}}(-1) 
\cdot u_{-\g}(1/x) u_{-\a{1}}(1)\\ 
&=& w^\prime {w^\prime}^{-1} w_{\a{1}}^{-1} {w^\prime} w_\g(x) 
\cdot n^\prime \cdot \bn \\
&=& w^\prime w_\g(d)^{-1} w_\g(x) \cdot n^\prime \cdot \bn \\
&=& w^\prime (w_\g \g^\vee(1/d))^{-1} w_\g \g^\vee(1/x) 
\cdot n^\prime \cdot \bn \\
&=& w^\prime \g^\vee(d/x) \cdot n^\prime \cdot \bn \\
&=& m n^\prime \bn.   
\end{eqnarray*} 

To see (\ref{m_alternate}), we use (\ref{conjugate_w_beta}) 
repeatedly to write 
\begin{eqnarray*} 
m &=& w_{\a{1}}^{-1} w^\prime w_\g(x) \\
&=& w_{\a{1}}^{-1} \cdot w^\prime w_\g(x) {w^\prime}^{-1} \cdot w^\prime \\ 
&=& w_{\a{1}}^{-1} w_{\a{1}}(D x) w^\prime \\
&=& w_{\a{1}}^{-1} w_{\a{1}} \a{1}^\vee(1/(D x)) w^\prime \\
&=& \a{1}^\vee(1/(D x)) w^\prime. 
\end{eqnarray*}
where 
\begin{equation} 
D =  
\begin{cases}   
\begin{array}{l}
d_{-\a{2},\b{2}} \cdots d_{-\a{n},\b{n}} d_{-\a{n+1},\b{n+1}} \\
\cdot \, d_{-\a{n},\b{n+2}} \cdots d_{-\a{2},\b{2n}}  
\end{array}
& \mbox{ if }\, \G \mbox{ is odd, } \\ 
 & \\ 
\begin{array}{l}
d_{-\a{2},\b{2}}  \cdots d_{-\a{n-1},\b{n-1}} d_{\a{n},\b{n}} \\ 
d_{-\a{n+1},\b{n+1}} d_{-\a{n-1},\b{n+2}} \cdots d_{-\a{2},\b{2n-1}} 
\end{array}
& \mbox{ if }\, \G \mbox{ is even. } \\ 
\end{cases} 
\end{equation}
Notice that conjugation by 
${w^\prime}^{-1}$ sends $\g=\b{2n}$ in the odd case and $\g=\b{2n-1}$ 
in the even case back to $\a{1}$.   

Finally, we claim that $D d = 1$ in both even and odd cases. To see 
this note that we can write 
\[ D d = 
\begin{cases}
\prod_{i=1}^n d_{\a{i+1},\b{i}} d_{-\a{i+1},\b{i+1}} 
\prod_{j=1}^{n-1} d_{\a{j+1},\b{2n-j}} d_{-\a{j+1},\b{2n+1-j}} 
& \mbox{ if }\, \G \mbox{ is odd, } \\ 
 & \\ 
\prod_{i=1}^n d_{\a{i+1},\b{i}} d_{-\a{i+1},\b{i+1}} 
\prod_{j=1}^{n-2} d_{\a{j+1},\b{2n-1-j}} d_{-\a{j+1},\b{2n-j}} 
& \mbox{ if }\, \G \mbox{ is even. } \\ 
\end{cases} 
\] 
Using (\ref{minus_a}) followed by (\ref{d_w_a}) we can rewrite this 
as 
\[ D d = 
\begin{cases}
\prod_{i=1}^n (-1)^{\langle \b{i}+\b{i+1}, \aa^\vee_{i+1}  \rangle}
\prod_{j=1}^{n-1} (-1)^{\langle \b{2n-j}+\b{2n+1-j}, \aa^\vee_{j+1}  \rangle} 
& \mbox{ if }\, \G \mbox{ is odd, } \\ 
 & \\ 
\prod_{i=1}^n (-1)^{\langle \b{i}+\b{i+1}, \aa^\vee_{i+1}  \rangle} 
\prod_{j=1}^{n-2} (-1)^{\langle \b{2n-1-j}+\b{2n-j}, \aa^\vee_{j+1} \rangle} 
& \mbox{ if }\, \G \mbox{ is even. } \\ 
\end{cases} 
\] 
Using the explicit root data we have described earlier we can 
see easily that in every single term of the above products the power 
of $(-1)$ is an even integer. In fact, $\b{i} + \b{i+1} \equiv \a{i+1} \mod 2$ 
for all $i$ and $\b{2n-j} + \b{2n+1-j} \equiv \a{j+1} \mod 2$ for all $j$ in 
the odd case and $\b{2n-1-j} + \b{2n-j} \equiv \a{j+1} \mod 2$ for all $j$ in 
the even case. This completes the proof. 
\end{proof} 
 
The following is assumption 4.1 of \cite{shahidi-imrn} for our cases. 
 
\begin{prop}\label{Assumption4.1} 
Let $n \in N$ satisfy (\ref{decomp}). Then except for a subset  
of measure zero of $N$ we have  
\begin{equation}\label{twosets}  
U_{M,n} = U^\prime_{M,m},  
\end{equation}  
where notation is as in Section 4 of \cite{shahidi-imrn}, i.e.,  
\begin{equation}  
U_{M,n} = \left\{ u\in U_M \,:\, u\, n\, u^{-1} = n\right\} 
\end{equation} 
and  
\[  
U^\prime_{M,m} = \left\{ u\in U_M \,:\, m\, u\, m^{-1} \in U_M  
\mbox{ and } \chi(m u m^{-1}) = \chi(m) \right\}. 
\] 
\end{prop} 
Note that the condition $\chi(m u m^{-1}) = \chi(m)$ in the  
definition of $U^\prime_{M,m}$ in our case is just the compatibility of $\chi$  
with elements of the Weyl group (cf. pages 2079--80 of \cite{shahidi-imrn}).

\begin{proof} 
By arguments such as those on page 2085 of \cite{shahidi-imrn} if  
the proposition is true for $n\in N$, then it is also true for  
every member of the intersection of its conjugacy class under  
$G$ with $N$ provided that for the $m$-part we take the twisted  
conjugacy classes (cf. (4.10) of \cite{shahidi-imrn}). Hence, it  
is enough to verify the proposition for those $n$ as in Corollary  
\ref{conjugacy}. Fix one such $n=u_{\a{1}}(1) u_\g(q)$ with  
$q \not= 0$ for the rest of this proof. 
 
We can explicitly compute both sides of (\ref{twosets}) as follows. 
Any $u \in U_M$ can be written as  
\begin{equation}\label{u_defn} 
u = \prod_{\bb \in \Sigma(\theta)^+} u_{\bb} (x_\bb),   
\end{equation} 
where the order of the terms in the product is with respect to the 
total order of $R$ we have fixed.  
 
We have $u n u^{-1} = n$ if and only if  
$u u_{\a{1}}(1) u_\g(q) u^{-1} = u_{\a{1}}(1) u_\g(q)$.  
Notice that by (\ref{comm}) we know that $u_\g(q)$ commutes with  
all $u_{\bb} (x_\bb)$ in the above product. Hence, $u\in U_{M,n}$ if  
and only if $u u_{\a{1}}(1) u^{-1} = u_{\a{1}}(1)$. Among  
the terms $u_{\bb} (x_\bb)$ the element $u_{\a{1}}(1)$ commutes  
with those with $\bb\in\Sigma(\Omega)^+$ where  
$\Omega = \Delta - \{\a{1},\a{2}\}$. Also, if  
$\bb\in\Sigma(\theta)^+ - \Sigma(\Omega)^+$, then so is  
$\g-\bb-\a{1}$. Using (\ref{comm}) several times we can now write  
\begin{eqnarray*} 
 u u_{\a{1}}(1) u^{-1}  
= \prod_{\bb \in \Sigma(\theta)^+} u_{\bb} (x_\bb) 
u_{\a{1}}(1) 
\left(\prod_{\bb \in \Sigma(\theta)^+} u_{\bb} (x_\bb)\right)^{-1} & \\  
= u_\g\left(\sum c_{\a{1},\bb;1,1} c_{\a{1}+\bb,\delta;1,1} x_\bb x_\delta \right) 
 \cdot \prod_{\bb \in \Sigma(\theta)^+ - \Sigma(\Omega)^+}  
u_{\a{1}+\beta}(- c_{\a{1},\beta} x_\bb) \cdot  
u_{\a{1}}(1),   
\end{eqnarray*} 
where the sum in the first term is over {\it unordered} pairs $(\bb,\delta)$ of  
roots in $\Sigma(\theta)^+ - \Sigma(\Omega)^+$ such that $\bb+\delta=\g-\aa$  
and $\bb\not=\delta$. Here the order of terms is prescribed by the order  
we fixed in (\ref{u_defn}). This implies that  
$u u_{\a{1}}(1) u^{-1} = u_{\a{1}}(1)$ if and only if  
$x_\bb = 0$ for all $\bb\in\Sigma(\theta)^+ - \Sigma(\Omega)^+$. Therefore,  
\begin{equation}\label{U_M_n}  
U_{M,n} = \left\{ \prod_{\bb \in \Sigma(\Omega)^+} u_{\bb} (x_\bb) \right\}.  
\end{equation}

To compute $U^\prime_{M,m}$ note that with $d$ as in (\ref{ddd}) we have   
\begin{eqnarray*} 
m\, u\, m^{-1} &=& w^\prime \g^\vee(d/q)  
                   \prod_{\bb \in \Sigma(\theta)^+} u_{\bb} (x_\bb) 
		   \g^\vee(d/q)^{-1} {w^\prime}^{-1} \\ 
&=& w^\prime \prod_{\bb \in \Sigma(\theta)^+} u_{\bb}  
(\bb(\g^\vee(d/q))\, x_\bb) {w^\prime}^{-1} \\  
&=& \prod_{\bb \in \Sigma(\theta)^+}  
w^\prime  u_{\bb} (\bb(\g^\vee(d/q))\, x_\bb) {w^\prime}^{-1}.  
\end{eqnarray*} 
 
Conjugation by the element $w^\prime$ sends each positive  
root group with root in $\Sigma(\theta)^+ - \Sigma(\Omega)^+$  
to a root group corresponding to a negative root and sends those with  
roots in $\Sigma(\Omega)^+$ to themselves. Therefore, again  
\begin{equation}\label{U_M,m} 
U^\prime_{M,m} = \left\{ \prod_{\bb \in \Sigma(\Omega)^+} u_{\bb} (x_\bb) \right\}. 
\end{equation} 
Now (\ref{twosets}) follows from (\ref{U_M_n}) and (\ref{U_M,m}). 
\end{proof} 
 
We would like to have an explicit identification of  
$\gl{1}(F)\times \wG{n}(F)$ with $M$ as a Levi subgroup of $G$. 
Going back to our descriptions of the groups  
$\gspin{2n+1}$ and $\wgspin{2n}$ in Section \ref{str},  
note that if we consider the root datum obtained from that of $\G$  
by eliminating $e_1$ and $e_1^*$ in the odd case and  
$E_1$ and $E_1^*$ in the even case as well as the root $\a{1}$ and its  
corresponding coroot, then it corresponds to a subgroup of  
$\G$ isomorphic to $\wG{n}$. Denote the $F$-points of this subgroup  
by $G^\sim_n$. Let $k\in G^\sim_n$ and let $a\in F^\times$. We claim that $e_1^*(a)$  
in the odd case and $E_1^*(a)$ in the even case and $k$ commute. To see  
this it is enough to observe that $e_1^*(a)$ or $E_1^*(a)$ 
commutes with $u_\bb(x)$ for all $\bb\in\Sigma(\theta)$ since $G^\sim_n$  
is generated by the corresponding $U_\bb$'s along with  
a subtorus of $T$. By (\ref{RootGroup}) we have  
\[ 
e_1^*(a) u_\bb(x) e_1^*(a)^{-1} = u_\bb(\bb(e_1^*(a)) x)  
= u_\bb(a^{\langle\bb,e_1^*\rangle} x)     
\] 
and similarly for $E_1^*(a)$. Moreover, $\langle\bb,e_1^*\rangle = 0$  
for all $\bb\in\Sigma(\theta)$. Therefore,  
$e_1^*(a)$ in the odd case and $E_1^*(a)$ in the even case and  
$u_\bb(x)$ commute. This implies that the map $(a,k)\mapsto e_1^*(a) k$  
in the odd case and $(a,k)\mapsto E_1^*(a) k$ in the even case  
is a homomorphism which gives the identification  
of $\gl{1}(F)\times \wG{n}(F)$ with $M$. In particular, the element  
$m = \a{1}^\vee(d/x) w^\prime$ in (\ref{m_alternate}) identifies with  
\begin{equation}\label{m-decomp} 
(d/x, e_2^*(x/d) w^\prime) \mbox{ or }  (d/x, E_2^*(x/d) w^\prime)  
\in \gl{1}(F)\times \wG{n}(F) 
\end{equation} 
since $\a{1}^\vee = e_1^* - e_2^*$ in the odd case and  
$\a{1}^\vee = E_1^* - E_2^*$ in the even case. (Notice that $w^\prime\in G^\sim_n$.)  
Moreover, $e_2^*(x/d) w^\prime$ or $E_2^*(x/d) w^\prime$ is an element of  
a maximal Levi subgroup in $\wG{n}(F)$ just as in the case of classical  
groups in \cite{shahidi-imrn}.  
 
We are now prepared to express the $\g$-factors as Mellin transforms. 
 
\begin{prop}\label{mellin} 
Let $\s$ be an irreducible admissible $\psi$-generic representation  
of $\wG{n}(F)$ (cf. Remark \ref{wgspin-center}).  
Consider $\gl{1}\times\wG{n}$ as a standard Levi in $\G$ as above.  
Let $\eta$ be any non-trivial character of $F^\times$ with $\eta^2$  
ramified. Then  
\begin{equation}\label{mellin-formula}  
\g(s,\eta\times\s,\psi)^{-1} = g(s,\eta) \cdot  
	 \int_{F^\times} j_{v,\overline{N}_0} (a(x)w^\prime)  
          \eta(x) |x|^{s-n-\delta} dx^\times,    
\end{equation}  
where $a(x)=e_2^*(x/d)$ or $E_2^*(x/d)$ are as in (\ref{m-decomp}),  
$\delta=1/2$ and $g(s,\eta)=\eta(-1)^n$ in the odd case and $\delta=1$  
and $g(s,\eta)=\eta(-1)^{n-1} \g(2s,\eta^2,\psi)^{-1}$ in the even case 
with $d$ as in (\ref{ddd}).  
Here $v\in V_\s$ and $W_v\in{\mathcal W}(\s,\psi)$ with $W_v(e)=1$,  
where ${\mathcal W}(\s,\psi)$ denotes the Whittaker model of $\s$.  
Moreover, $\overline{N}_0\subset\overline{N}$ is a sufficiently large  
compact open subgroup of the opposite unipotent subgroup $\overline{N}$  
to $N$, where $P=MN\subset G$ is the Levi decomposition of  
the corresponding standard parabolic subgroup. The function  
$j_{v,\overline{N}_0}$ denotes the partial Bessel function  
defined in \cite{shahidi-imrn}.  
\end{prop} 
 
\begin{proof}  
Given that $\eta^2$ is ramified this proposition is the main result  
of \cite{shahidi-imrn}, (6.39) of Theorem 6.2, applied to our cases. Notice  
that the two hypotheses of that theorem, i.e., Assumptions 4.1 and 5.1  
of \cite{shahidi-imrn}, for our cases are our Propositions  
\ref{Assumption4.1} and \ref{Assumption5.1}, respectively. 
 
To get from (6.39) of \cite{shahidi-imrn} to (\ref{mellin-formula})  
above note that we have  
\[ \omega_{\s_s}^{-1}(\dot{x}_\aa) (w_0\omega_{\s_s})(\dot{x}_\aa)  
= \eta(x)^2 |x|^{2s} \] 
Moreover, as in Section 7 of \cite{shahidi-imrn} we have  
\[ q^{\langle s \tilde{\aa}+\rho,H_M(\dot{m})\rangle} d\dot{n} =  
|x|^{-s-n+\delta} dx^\times \] and  
\[ j_{\tilde{v},\overline{N}_0}(\dot{m}) = \eta(d/x)  
j_{v,\overline{N}_0} (a(x)w^\prime), \] 
where $a(x)=e_2^*(x/d)$ in the odd case and $E_2^*(x/d)$ in the even case. 
\end{proof} 
 
In the next section we first rewrite (\ref{mellin-formula})  
in terms of Bessel functions defined as in  \cite{cogdell-ps} and  
then study their asymptotics. 
 
\subsection{Bessel functions and their asymptotics}\label{bess-sec} 
 
We now briefly review some basic facts from \cite{cogdell-ps}. In  
view of \cite{ckpss-classical} and  
particularly \cite{stability1,stability2} which will study these issues in  
more generality, we only concentrate on the cases at hand in this article  
and leave out the details of the more general situation, thereby simplifying  
some of the notations. 
 
We will use the same notation as \cite{cogdell-ps}. Consider  
the group $\wG{n}$ in both even and odd cases and $w^\prime$ above as  
an element of its Weyl group. Notice that $w^\prime$ supports a Bessel  
function and is minimal (in Bruhat order) non-trivial with respect  
to this property. Moreover,  
$A_{w^\prime} = Z_{M_\Omega}$, where $\Omega = \Delta - \{\a{1},\a{2}\}$,  
$M_\Omega$ is the standard Levi determined by $\Omega$, and $Z_{M_\Omega}$  
denotes its center.

Let $\s$ be an irreducible admissible $\psi$-generic representation of  
$\wG{n}(F)$ and take $v\in V_\s$ such that the associated  
Whittaker function $W_v\in{\mathcal W}(\s,\psi)$ satisfies $W_v(e)=1$.  
The associated Bessel function on $Z_{M_\Omega}$ is defined via  
\begin{equation} 
J_{\s,w^\prime}(a) = \int_{U^-_{w^\prime}} W_v(a w^\prime u) \psi^{-1}(u) du,  
\end{equation}  
where $a\in Z_{M_\Omega}$ and $U_{w^\prime}^{-} = \prod_{\aa} \,  U_\aa$,  
where the product is over all those $\aa\in\Sigma(\theta)^+$ with $w^\prime(\aa) < 0$  
and $U_\aa$ is as before. 
 
Similar to \cite{cogdell-ps,ckpss-classical} we have that $J_{\s,w^\prime}$  
exists and is independent of $v\in V_\s$ and for convergence purposes we use  
a slight modification of it, namely, the partial Bessel function  
\begin{equation} 
J_{\s,w^\prime,v,Y} (a) = \int_Y W_v(a w^\prime y) \psi^{-1}(y) dy,  
\end{equation}  
where $Y\subset U_{w^\prime}^{-}$ is a compact open subgroup. 
 
\subsubsection{Domain of Integration} 
 
We now show that the partial Bessel functions of \cite{shahidi-imrn} are  
the same as those in \cite{cogdell-ps}. 
 
Recall that $M=M_\theta=\gl{1}(F)\times\wG{n}(F)\subset G=\G(F)$  
and consider $m \in M$ as in (\ref{m-decomp}), i.e., $m=(d/x,e_2^*(x/d) w^\prime)$  
in the odd case and $m=(d/x,E_2^*(x/d) w^\prime)$ in the even case with 
$d$ as in (\ref{ddd}).  
We now prove that the Bessel functions of \cite{shahidi-imrn} and those  
of \cite{cogdell-ps} are actually the same. 
 
\begin{lem}  
We can choose $Y$ appropriately so that with  
$m^\prime = e_2^*(x/d)w^\prime$ or $E_2^*(x/d)w^\prime$ as above we have   
\[ j_{v,\bN_0} (m^\prime) =  
\begin{cases} 
J_{\s,w^\prime,v,Y} (e_2^*(x/d))  
& \mbox{ in the odd case } \\  
J_{\s,w^\prime,v,Y} (E_2^*(x/d))  
& \mbox{ in the even case. }  
\end{cases} \]  
Here, $m^\prime$ is an element of a maximal Levi  
$M^\prime = \gl{1}(F) \times \wG{n-1}(F)$ in $\wG{n}(F)$.  
\end{lem} 
 
\begin{proof} 
Let us first recall Theorem 6.2 of \cite{shahidi-imrn}. In the notation  
of that paper we have $j_{v,\bN_0}(m^\prime) = j_{v,\bN_0}(m^\prime,y_0)$ with  
$y_0\in F^\times$ satisfying  
$\mbox{ord}_F(y_0) = - \mbox{cond}(\psi) - \mbox{cond}(\eta^2)$. Here,  
the function $j_{v,\bN_0}(m^\prime,y_0)$ is given by  
\begin{equation}\label{bessel} 
\int_{U_{M^\prime,n}\backslash U_{M^\prime}} W_v(m^\prime u^{-1})  
\phi(u \aa^\vee(y_0)^{-1} \aa^\vee(x_\aa) \bn \aa^\vee(x_\aa)^{-1}  
\aa^\vee(y_0) u^{-1})  
\psi(u) du,  
\end{equation} 
where $\phi$ is the characteristic function of $\bN_0$,  
$x_\aa = 1/x$, and $\bn$ is as in Proposition \ref{conj-reps}.  
Again as in \cite{cogdell-ps,ckpss-classical} it follows from  
Proposition \ref{Assumption5.1} that we can take $U_{M,n}\backslash U_M$  
to be $U^-_{w^\prime}$. Notice that this only depends on $w^\prime$. On  
the other hand, $u\in U^-_{w^\prime}$ is in the domain of integration  
if and only if  
\[ u \aa^\vee(y_0)^{-1} \aa^\vee(x_\aa) \bn  
\aa^\vee(x_\aa)^{-1} \aa^\vee(y_0) u^{-1} \in \bN_0. \]  
This condition is equivalent to  
$u \aa^\vee(x_\aa) \bn \aa^\vee(x_\aa)^{-1} u^{-1}  
\in \aa^\vee(y_0)\bN_0\aa^\vee(y_0)^{-1}$ and  
$\aa^\vee(y_0)\bN_0\aa^\vee(y_0)^{-1}$ is another compact  
open subgroup of the same type as $\bN_0$ which we may replace  
it with.  
 
Recall that an arbitrary element of $\bN$ is given by  
\begin{equation} 
\bn(y)=\prod_{\aa\in R(\bN)} u_\aa(y_\aa),  
\end{equation} 
where $y = (y_\aa)_{\aa\in R(\bN)}$ and $y_\aa\in F$.  
Also recall that $\bn$ in (\ref{bessel})  
was given by $\bn = u_{-\g}(1/x) u_{-\a{1}}(1)$.  
Moreover, note that $x_\aa = 1/x$.  
Hence, $\aa^\vee(x_\aa) \bn \aa^\vee(x_\aa)^{-1} = \bn(y^\prime)$,  
where $y^\prime_{-\g} = 1$ and $y^\prime_{-\a{1}} = x$ and  
all other coordinates of $y^\prime$ are zero (cf. (\ref{RootGroup})).   
Of course, throughout we have a fixed ordering of the roots in the  
products similar to that of (\ref{n-expr}). Hence, the domain of  
integration is determined by $u \bn(y^\prime) u^{-1} \in \bN_0$. 
 
We may take $\bN_0 = \{\bn(y) : y_\aa\in\mathfrak p^{M_\aa}\}$ for  
all $\aa\in R(\bN)$ for sufficiently large integer vector  
$M =(M_\aa)_{\aa\in R(\bN)}$. As the $M_\aa$ increase, $\bN_0$ will  
exhaust $\bN$.  
 
On the other hand, any $u\in U^-_{w^\prime}$ is given by  
\begin{equation} 
u = u(b) = \prod_{\aa\in\Sigma^+(\Omega)} u_\aa(b_\aa) 
\end{equation} 
for $b = (b_\aa)_{\aa\in\Sigma^+(\Omega)}$ and $b_\aa\in F$.  
Now, $u \bn(y^\prime) u^{-1} = \bn(y^{\prime\prime})$  
where $y^{\prime\prime}$ depends linearly  
on $b$ and $y^\prime$. In other words, $y^{\prime\prime}$ depends  
upon $x$ and $b$. Of course, we could compute $y^{\prime\prime}$  
explicitly in terms of $x$ and $b$ using structure constants; however,  
that will have no bearing on what follows and is not needed. Now,  
choose $Y=\{u=u(b) : y^{\prime\prime} \geq M\}$. This defines the  
domain of integration. Enlarging $\bN_0$ if need be will then  
imply that $j_{v,\bN_0} (m)$ does not depend on $m$ and we  
conclude the lemma. 
\end{proof} 
 
Therefore, we can rewrite (\ref{mellin-formula}) as  
\begin{equation}\label{mellin-formula2}  
\g(s,\eta\times\s,\psi)^{-1} = g(s,\eta) \cdot  
	 \int_{F^\times}  
J_{\s,w^\prime,v,Y} (a(x)) \eta(x) |x|^{s-n-\delta} dx^\times,  
\end{equation}  
with $a(x) = e_2^*(x/d)$ or $E_2^*(x/d)$. 
 
\subsubsection{Asymptotics of Bessel functions}  
 
We now study the asymptotics of our Bessel functions near zero  
and infinity. This will allow us to prove our result on  
stability given that the $\g$-factors are already written  
as the Mellin transform of Bessel functions.  
 
Our starting point is the following analogue of proposition 5.1  
of \cite{cogdell-ps} for our groups. Note that the proposition  
was only proved for the group $\so{2n+1}$. However, as was  
pointed out in \cite{ckpss-classical}, the methods used to  
prove it are quite general. This was pointed out for classical  
groups (with finite center) in \cite{ckpss-classical} but  
the same also holds for $\gspin{}$ or $\wgspin{}$ groups which  
are of interest  
to us since the only difference is the infinite center which  
is already contained in the fixed Borel subgroup we are modding  
out with. 
 
\begin{prop}\label{asymp1} 
There exists a vector $v^\prime_\s\in V_\s$ and a compact neighborhood  
$BK_1$ of the identity in $B\backslash G^\sim_n$ such that if $\chi_1$  
is the characteristic function of $BK_1$, then  
for all sufficiently large compact open sets $Y\subset U_w^-$ we  
have  
\begin{equation} 
J_{\s,w,v,Y} (a) = \int_Y W_v(a w y) \chi_1(a w y) \psi^{-1}(y) dy +  
W_{v^\prime_\s}(a).  
\end{equation} 
 
Notice that $w$ here would be our earlier $w^\prime$ if we want to  
consider the group $G^\sim_n$ as part of the Levi $M$ in $G$  
as in (\ref{m-decomp}). 
\end{prop} 
 
We now would like to rewrite this in a way that only depends on the  
central character of $\s$. To this end we argue similar to \cite{ckpss-classical}  
making some necessary modifications along the way. For any positive  
integer $M$ set  
\[ U(M) = \langle u_\aa(t) : \aa\in\Delta, |t| \leq q^M \rangle. \]  
These are compact open subgroups of $U$  and as the integer $M$ grows, they  
exhaust $U$. For any $v\in V_\s$ we define  
\[ v_M = \frac{1}{\mbox{Vol}(U(M))} \int_{U(M)} \psi^{-1}(u) \s(u)v du. \]  
Smoothness of $\s$ implies that this is a finite sum and  
$v_M\in V_\s$. Then just as in \cite{cogdell-ps} if $Y$ is sufficiently  
large relative to $M$, then we may choose $v^\prime_\s$  
and $K_1$ such that $K_1\subset\mbox{Stab}(v_M)$ and we have  
\begin{equation} 
\int_Y W_v(a w y)\chi_1(a w y)\psi^{-1}(y) dy = \int_Y W_{v_M}(a w y)  
\chi_1(a w y)\psi^{-1}(y) dy. 
\end{equation} 
 
Write $a w y = u t k_1$ with $u\in U$, $t\in T$ and $k_1\in K_1$. Then  
$K_1\subset\mbox{Stab}(v_M)$ implies that $W_{v_M}(a w y)=\psi(u) W_{v_M}(t)$.  
As in \cite{cogdell-ps} the support of $W_{v_M}$ is contained in  
\begin{equation} 
T_M = \left\{ t \in T : \aa(t) \in 1 + {\mathfrak p}^M  
\mbox{ for all simple }\aa \right\}. 
\end{equation} 
 
At this point we would like to assume that the center $\Z$ is connected which  
we can do (cf. Remark \ref{wgspin-center}). This is why we chose to work with  
the group $\wgspin{}$ in the even case in this section. By connectedness of $\Z$  
and since the groups are split we have the following exact sequence of $F$-points  
of tori 
\begin{equation} 
0 \vlto{30} Z \vlto{30} T \vlto{30} T_{ad} \vlto{30} 0 
\end{equation} 
which splits (cf. \cite{stability1,stability2}). Recall that $Z=\Z(F)$ and so on.  
Identify $T_{ad}$ with $(F^\times)^n$ through values of roots  
and let $T_M^1 \subset T$ be the image of  
\[ (1+{\mathfrak p}^M)^n \subset  T_{ad}, \] 
under the splitting map. Here, the rank of $T_{ad}$ is $n$ and $T_M = Z T_M^1$.  
 
Now if $t\in T$, then we can write $t = z t^1$ with $z\in Z$ and  
$t \in T_M^1$. Also we have $W_v(t) = W_{v_M}(t)$ and if we choose $M$  
large enough so that $T_M^1\subset T\cap\mbox{Stab}(v)$, then  
\[ W_{v_M}(t) = W_v(t) = W_v (z t^1) = \omega_\s(z) W_v(t^1) = \omega_\s(z). \]  
 
Therefore, in our integral, $W_{v_M}(a w y)\chi_1(a w y) \not= 0$ if and only if  
$a w y \in U T_M K_1$ or $y\in (a w)^{-1} U T_M K_1$. Writing  
$a w y = u t k_1 = u(a w y) z(a w y) t^1 k_1$ then implies that  
\begin{equation} 
\underset{Y}{\int} W_v(a w y)\chi_1(a w y)\psi^{-1}(y) dy =  
\underset{Y\cap (a w)^{-1}U T_M K_1}{\int} \psi(u(a w y))  
\psi^{-1}(y) \omega_\s(z(a w y)) dy.   
\end{equation} 
 
Therefore, we can rewrite Proposition \ref{asymp1} as follows. 
 
\begin{prop}\label{asymp} 
Let $v\in V_\s$ with $W_v(e) = 1$ and choose $M$ sufficiently large  
so that $T_M^1 \subset T \cap \mbox{Stab}(v)$. There exists a vector  
$v_\s^\prime \in V_\s$ and a compact open subgroup $K_1$ such that  
for $Y$ sufficiently large we have  
\[ J_{\s,w,v,Y} (a) = \int_{Y\cap (aw)^{-1}U T_M K_1}  
\psi(u(a w y)) \psi^{-1}(y) \omega_\s(z(a w y)) dy +  
W_{v^\prime_\s}(a). \] 
\end{prop}

\subsection{Proof of Theorem \ref{stable}} 
We are now prepared to prove Theorem \ref{stable}.  
\begin{proof} 
Let $\s_i=\pi_i$, $i=1,2$ in the odd case. In the even case choose  
a character $\mu$ of the center $Z^\sim$ of $\wgspin{2n}(F)$ (which  
contains the center of $\gspin{2n}(F)$) such that $\mu$ agrees with the  
central characters $\oo_{\pi_1}=\oo_{\pi_2}$ on the center of  
$\gspin{2n}(F)$. Consider the representation of $\wgspin{2n}(F)$  
induced from the representation $\mu\otimes\pi_i$ on $Z^\sim \cdot \gspin{2n}(F)$  
(which is of finite index in $\wgspin{2n}(F)$) and let $\s_i$ be an  
irreducible constituent of this induced representation (cf. \cite{tadic}).  
Note that the choice of $\s_i$ is irrelevant.  
Then, \[ \g(s,\eta\times\s_i,\psi) = \g(s,\eta\times\pi_i,\psi), \] 
by Remark \ref{wgspin-center}. Also,  
the assumption $\oo_{\pi_1} = \oo_{\pi_2}$ implies  
$\oo_{\s_1} = \oo_{\s_2}$.

Choose $v_i\in V_{\sigma_i}$, $i=1,2$ with $W_{v_i} (e) =1$ and  
let $M$ be a large enough integer so that  
$T_M^1 \subset T \cap \mbox{Stab}(v_i)$. Choose  
a compact open subgroup $K_0 \subset \mbox{Stab}(v_1) \cap \mbox{Stab}(v_2)$.  
Then in Proposition \ref{asymp} we may take  
\[ K_1 = \bigcap_{u\in U(M)} u^{-1} K_0 u, \] 
i.e., we can take the same $K_1$ for both $\sigma_1$ and $\sigma_2$. Consequently, by  
Proposition \ref{asymp} there exist $v_{\sigma_i}^\prime \in V_{\sigma_i}$ such  
that  
\begin{equation} 
J_{\sigma_i,w,v,Y} (a) = \int_{Y\cap (aw)^{-1}U T_M K_1}  
\psi(u(awy)) \psi^{-1}(y) \omega_{\sigma_i}(z(awy)) dy +  
W_{v^\prime_{\sigma_i}}(a). 
\end{equation} 
Now $\omega_{\sigma_1} = \omega_{\sigma_2}$ implies that  
\begin{equation} 
J_{\sigma_1,w,v,Y} (a) - J_{\sigma_2,w,v,Y} (a) =  
W_{v^\prime_{\sigma_1}}(a) - W_{v^\prime_{\sigma_2}}(a). 
\end{equation}

Now taking $a=a(x)$ to be $e_2^*(x/d)$ or $E_2^*(x/d)$ and  
$w$ to be $w^\prime$ described before, we  
apply (\ref{mellin-formula2}) to conclude that  
\begin{eqnarray*}  
& \g(s,\eta\times\s_1,\psi)^{-1} - \g(s,\eta\times\s_2,\psi)^{-1} \\ 
& \\ 
&=g(s,\eta) \int_{F^\times} (J_{\s_1,w,v,Y} (a(x)) - J_{\s_2,w,v,Y} (a(x)) )  
\,\, \eta(x) \,\, |x|^{s-n+\delta} d^\times x \\ 
& \\ 
&=g(s,\eta) \int_{F^\times} ( W_{v^\prime_{\s_1}} (a(x)) - W_{v^\prime_{\s_2}} (a(x)))  
\,\, \eta(x) \,\, |x|^{s-n+\delta} d^\times x.   
\end{eqnarray*} 
 
However, note that Whittaker functions are smooth and for $\Re(s) >> 0$ and $\eta$  
sufficiently ramified we have  
\[ 
\int_{F^\times} W_{v^\prime_{\s_i}} (a(x)) \eta(x)  
|x|^{s-n+\delta} d^\times x \equiv 0. 
\] 
 
Hence, for $\Re(s) >> 0$ we have  
$\g(s,\eta\times\s_1,\psi)^{-1} - \g(s,\eta\times\s_2,\psi)^{-1} \equiv 0$ which  
then implies $\g(s,\eta\times\s_1,\psi)=\g(s,\eta\times\s_2,\psi)$ for all $s$ by  
analytic continuation. Therefore, $\g(s,\eta\times\pi_1,\psi)=\g(s,\eta\times\pi_2,\psi)$. 
\end{proof}

\subsection{Stable Form of $\g(s,\eta\times\pi,\psi)$} 
 
We now prove some consequences of Theorem \ref{stable} which are important to us later. 
 
First, let us compute the stable form of Theorem \ref{stable} by taking  
$\pi_2$ to be an appropriate principal series representation and computing  
its right hand side explicitly. 
 
\begin{prop}\label{StableGamma} 
Let $\pi$ be an irreducible generic representation of $\G_n(F)$ with central character  
$\oo=\oo_\pi$. Let $\mu_1,\dots,\mu_n$ be $n$ character of $F^\times$.  
Then for every sufficiently ramified character $\eta$ of $F^\times$ we have  
\[ 
\g(s,\eta\times\pi,\psi) = \prod_{i=1}^n \g(s,\eta\mu_i,\psi) 
                                               \g(s,\eta\oo\mu_i^{-1},\psi). 
\] 
\end{prop} 
\begin{proof} 
Set $\mu_0=\oo$ and consider the character  
\[ \mu =(\mu_0\circ e_0) \otimes (\mu_1\circ e_1) \otimes \cdots  
\otimes (\mu_n\circ e_n) \]  
of $\T(F)$ with $e_i$'s as in Section \ref{structure}. Proposition  
\ref{center-gspin} implies that the restriction of the character  
$\mu$ to the center of $\G_n(F)$ is $\mu_0=\oo$. Consider the induced  
representation $\mbox{Ind}(\mu)$ from the Borel to $\G_n(F)$. Reordering the $\mu_i$ 
if necessary, we may assume that it has an irreducible admissible  
generic subrepresentation $\pi_2$ (cf. Proposition \ref{smc}).  
Since $\oo_{\pi_2}=\mu_0=\oo=\oo_\pi$, we  
can apply Theorem \ref{stable} to get  
$\g(s,\eta\times\pi,\psi) = \g(s,\eta\times\pi_2,\psi)$. Multiplicativity  
of $\g$-factors can now be used to compute the right hand side to get  
\begin{equation}\label{mult}  
\g(s,\eta\times\pi_2,\psi) =  
\prod_{i=1}^n \g(s,\eta\mu_i,\psi)\,\, \g(s,\eta\oo\mu_i^{-1},\psi)  
\end{equation} 
which finishes the proof.  
\end{proof} 
\begin{cor}\label{StableFormL-Ep} 
Let $\pi$ be an irreducible generic representation of $\G_n(F)$ with central character  
$\oo=\oo_\pi$. Let $\mu_1,\dots,\mu_n$ be $n$ character of $F^\times$ as in  
Proposition \ref{StableGamma}. Then for every sufficiently ramified  
character $\eta$ of $F^\times$ we have  
\[ L(s,\eta\times\pi) \equiv 1 \] and 
\[ 
\e(s,\eta\times\pi,\psi) = \prod_{i=1}^n \e(s,\eta\mu_i,\psi)\e(s,\eta\oo\mu_i^{-1},\psi). 
\] 
\end{cor} 
\begin{proof} 
If $\eta$ is sufficiently ramified, then by \cite{shahidi-twist} we have  
\[ L(s,\eta\times\pi) \equiv 1. \] 
This implies that $\e(s,\eta\times\pi,\psi)=\g(s,\eta\times\pi,\psi)$. Moreover,  
since $\eta$ is highly ramified so is each $\eta\mu_i$ and $\eta\oo\mu_i^{-1}$  
which implies that $L(s,\eta\mu_i)\equiv 1$ and $L(s,\eta\oo\mu_i^{-1})\equiv 1$.  
Therefore, $\e(s,\eta\mu_i,\psi)=\g(s,\eta\mu_i,\psi)$ and  
$\e(s,\eta\oo\mu_i^{-1},\psi)=\g(s,\eta\oo\mu_i^{-1},\psi)$. Now the second statement  
of the corollary follows from Proposition \ref{StableGamma}. 
\end{proof} 
 
\section{Analytic properties of global $L$-functions}\label{global}  
 
In this section we prove the properties of global $L$-functions that we need  
in order to apply the converse theorem. 
 
We again let $\G_n$ denote either the group $\gspin{2n+1}$ or $\gspin{2n}$ as  
in Section \ref{structure}. Let $k$ be a number field and let $\Ad$ denote its  
ring of adeles. Let $S$ be a finite set of finite places of $k$. Let $\cT(S)$  
denote the set of irreducible cuspidal automorphic representations $\tau$ of  
$\gl{r}(\Ad)$ for $1\leq r \leq N-1$ such that $\tau_v$ is unramified for all  
$v\in S$. If $\eta$ is a continuous complex character of $k^\times\backslash\Ad^\times$,  
then let  
$\cT(S;\eta)=\left\{\tau=\tau^\prime\otimes\eta\,:\,\tau^\prime\in\cT(S)\right\}$.  
 
If $\pi$ is a globally generic cuspidal representation of $\G_n(\Ad)$ and $\tau$  
is a cuspidal representation of $\gl{r}(\Ad)$ in $\cT(S;\eta)$, then  
$\s = \tau\otimes\tilde{\pi}$ is a (unitary) cuspidal globally generic representation  
of $\M(\Ad)$, where $\M = \gl{r} \times \G_n$ is a Levi subgroup of a standard  
parabolic subgroup in $\G_{r+n}$. The machinery of the Langlands-Shahidi method  
as in Section \ref{local} now applies \cite{shahidi:88annals, shahidi:90annals}.  
Recall that  
\begin{eqnarray} 
L(s,\pi\times\tau) = \prod_{v} L(s,\pi_v\times\tau_v), \\ 
\e(s,\pi\times\tau) = \prod_{v} \e(s,\pi_v\times\tau_v, \psi_v),  
\end{eqnarray} 
where the local factors are as in (\ref{RS-L}) and (\ref{RS-ep}). 
 
\begin{prop}\label{entire} 
Let $S$ to be a non-empty set of finite places of $k$ and let $\eta$ be a  
character of $k^\times\backslash\Ad^\times$ such that $\eta_v$ is highly  
ramified for $v\in S$. Then for all $\tau\in\cT(S;\eta)$ the L-function  
$L(s,\pi\times\tau)$ is entire.  
\end{prop} 
\begin{proof} 
These $L$-functions are defined via the Langlands-Shahidi method as we outlined  
in Section \ref{local}. Now, the proposition is a special case of a more general  
result, Theorem 2.1 of \cite{kim-sha}. Note that we have proved the necessary  
assumption of that theorem, Assumption 1.1 of \cite{kim-sha}, for our cases in  
Proposition \ref{assumptionA}.   
\end{proof} 
 
The following lemma is an immediate consequence of Proposition \ref{assumptionA}.  
\begin{lem}\label{globalAssumptionA}  
The global normalized intertwining operator $N(s,\s,w)$ is a holomorphic and  
non-zero operator for $\Re(s) \geq 1/2$.  
\end{lem} 
 
\begin{prop}\label{bdd} For any cuspidal automorphic representation $\tau$  
of $\gl{r}(\Ad_F)$, $1\leq r \le 2n-1$, the $L$-function $L(s,\pi\times\tau)$  
is bounded in vertical strips.  
\end{prop} 
\begin{proof} 
This follows as a consequence of Theorem 4.1 of \cite{gel-sha} along the  
lines of Corollary 4.5 thereof, given the fact that we have proved  
Assumption 2.1 of \cite{gel-sha} in our Proposition \ref{globalAssumptionA}  
for our cases.  
\end{proof} 
\begin{prop}\label{fneq} For any cuspidal automorphic representation $\tau$  
of $\gl{r}(\Ad_F)$, $1\leq r \le 2n$, we have the functional equation  
\[ L(s,\pi\times\tau) = \epsilon(s,\pi\times\tau) L(1-s,\w{\pi}\times\w{\tau}). \]    
\end{prop} 
\begin{proof} 
This is a special case of Theorem 7.7 of \cite{shahidi:90annals}. 
\end{proof} 
 
\section{Proof of Main Theorem}\label{proof} 
 
As mentioned before, we will use the following variant of converse  
theorems of Cogdell and Piatetski-Shapiro. 
 
\begin{thm}\label{conversethm}  
Let $\Pi=\otimes^\prime\Pi_v$ be an irreducible admissible representation  
of $\gl{N}(\Ad)$ whose central character $\oo_\Pi$ is invariant under  
$k^\times$ and whose $L$-function $L(s,\Pi)=\prod_{v} L(s,\Pi_v)$ is  
absolutely convergent in some right half plane. With notation as in  
Section \ref{global}, suppose that for every $\tau\in\cT(S;\eta)$ we have 
\begin{itemize} 
\item[(1)] $L(s,\Pi\times\tau)$ and $L(s,\widetilde{\Pi}\times\widetilde{\tau})$  
extend to entire functions of $s\in\cx$, 
\item[(2)] $L(s,\Pi\times\tau)$ and $L(s,\widetilde{\Pi}\times\widetilde{\tau})$  
are bounded in vertical strips, and 
\item[(3)] $L(s,\Pi\times\tau)=\epsilon(s,\Pi\times\tau)L(1-s,\widetilde{\Pi} 
\times\widetilde{\tau})$. 
\end{itemize} 
Then, there exists an automorphic representation $\Pi^\prime$ of $\gl{N}(\Ad)$  
such that $\Pi_v\simeq\Pi_v^\prime$ for all $v\not\in S$. 
 
Here, the twisted $L$- and $\epsilon$-factors are defined via  
\[ L(s,\Pi\times\tau)=\prod_{v} L(s,\Pi_v\times\tau_v) \quad  
\epsilon(s,\Pi\times\tau)=\prod_{v}\epsilon(s,\Pi_v\times\tau_v,\psi_v)\] 
with local factors as in \cite{converse1}. 
\end{thm} 
 
This is the exact variant of converse theorems that appeared in  
Section 2 of \cite{cogdell-kim-ps-sha}.  
 
We can now prove Theorem \ref{transfer}. 
 
\begin{proof} --- We apply Theorem \ref{conversethm} with $N=2n$.  
We continue to denote by $\G_n$ either of $\gspin{2n+1}$ or $\gspin{2n}$. 
First, we introduce a candidate for the representation $\Pi$. Consider  
$\pi=\otimes^\prime \pi_v$ and let $S$ be as in the statement of the  
theorem, i.e., a non-empty set of non-archimedean places $v$ such that  
for all finite $v\not\in S$ both $\pi_v$ and $\psi_v$ are unramified. 
 
\begin{itemize} 
\item[{\bf(i)}] $v<\infty$ and $\pi_v$ unramified: Choose $\Pi_v$  
as in the statement of the theorem via the Frobenius-Hecke  
(or Satake) parameter.  
 
More precisely, since $\pi_v$ is unramified, it is given by an unramified  
character $\chi$ of the maximal torus $\T(k_v)$.  
This means that there are unramified characters  
$\chi_0,\chi_1,\dots,\chi_n$ of $k_v^\times$ such that for $t\in\T(k_v)$  
\begin{equation}\label{UnramifiedChar} 
 \chi(t)=(\chi_0\circ e_0)(t) (\chi_1\circ e_1)(t) \cdots (\chi_n\circ e_n)(t), 
\end{equation} 
where $e_i$'s form the basis of the rational characters of the maximal  
torus of $\G$ as in Section \ref{structure}. The character $\chi$  
corresponds to an element $\hat{t}$ in $\hT$, the maximal torus of  
(the connected component) of the Langlands dual group which is  
$\gsp{2n}(\cx)$ or $\gso{2n}(\cx)$, uniquely determined by the equation  
\begin{equation}\label{satake} 
\chi(\phi(\varpi)) = \phi(\hat{t}),  
\end{equation}  
where $\varpi$ is a uniformizer of our local field $k_v$ and  
$\phi\in X_*(\T)=X^*(\hT)$ (cf.(I.2.3.3) on page 26 of \cite{gelbart-sha-book}).  
We make this identification explicit via the correspondence  
$e_i^* \longleftrightarrow e_i$ for $i=0,\dots,n$  as in Section  
\ref{structure} which gave the duality of $\gspin{2n+1}\longleftrightarrow\gsp{2n}$  
and $\gspin{2n}\longleftrightarrow\gso{2n}$. Applying (\ref{satake})  
with the $\phi$ on the left replaced with $e_i^*$ and the one on the right  
replaced with $e_i$ for $i=0,1,\cdots, n$ yields  
\begin{equation}\label{values} 
\chi_i(\varpi) = \chi(e_i^*(\varpi)) = e_i(\hat{t}), \quad i=0,1,\cdots, n. 
\end{equation} 
 
We can now compute the Satake parameter explicitly as an element $\hat{t}$  
in the maximal torus $\hT$ of $\gsp{2n}(\cx)$ or $\gso{2n}(\cx)$, as  
described in (\ref{torus-gsp-gso}). If we write our unramified characters  
as $\chi_i(\ ) = |\ |_v^{s_i}$ for $\s_i\in\cx$ and $0 \le i \le n$, then we get  
\begin{equation}\label{parameter}  
	\hat{t}=\left(    \begin{array}{cccccc}  
                                   |\varpi|^{s_1}&   & & & & \\ 
                                      &\ddots&&&& \\  
                                   &&|\varpi|^{s_n}&&& \\ 
                                   &&&|\varpi|^{s_0-s_n}&& \\ 
                                   &&&&\ddots& \\ 
                                   &&&&&|\varpi|^{s_0-s_1}  
                      \end{array}   \right). 
\end{equation}  
Hence, $\Pi_v$ is the unique unramified constituent of the representation of  
$\gl{2n}(k_v)$ induced from the character  
\begin{equation}\label{RepPi_v} 
\chi_1\otimes\cdots\otimes\chi_n\otimes\chi_0\chi_n^{-1} 
\otimes\cdots\otimes\chi_0\chi_1^{-1}  
\end{equation} 
of the $k_v$-points of the standard maximal torus in $\gl{2n}$.  
 
A crucial point here is what the central characters of $\pi_v$ and $\Pi_v$  
are. It follows from Proposition \ref{center-gspin} that the central character  
$\oo_{\pi_v}=\chi_0$. Moreover, the central character $\oo_{\Pi_v}$ of $\Pi_v$  
is $\chi_0^n$, hence we have $\oo_{\Pi_v} = \oo_{\pi_v}^n$. 
 
Furthermore, note that $\w{\Pi}_v$ is the  
unique unramified constituent of the representation induced from  
\[ \chi_1^{-1}\otimes\cdots\otimes\chi_n^{-1}\otimes 
\chi_0^{-1}\chi_n\otimes\cdots\otimes\chi_0^{-1}\chi_1. \] 
Therefore, we have $\w{\Pi}_v\simeq\chi_0^{-1}\otimes\Pi_v$. In other  
words, $\Pi_v\simeq\w{\Pi}_v\otimes\omega_{\pi_v}$.  
 
\item[{\bf(ii)}] $v|\infty$: Choose $\Pi_v$ as in the statement of  
Theorem \ref{transfer} \cite{langlands-real}. 
 
To be more precise, Langlands associates to $\pi_v$ a homomorphism  
$\phi_v$ from the local Weil group $W_v=W_{k_v}$ to the dual group $\hG$  
which is $\gsp{2n}(\cx)$ or $\gso{2n}(\cx)$ in our cases. Both of  
these groups have natural embeddings $\iota$ into $\gl{2n}(\cx)$ and  
we take $\Pi_v$ to be the irreducible admissible representation  
of $\gl{2n}(k_v)$ associated to $\Phi_v = \phi_v\circ\iota$. 
 
We want to show that again we have $\oo_{\Pi_v} = \oo_{\pi_v}^n$ and  
$\Pi_v\simeq\w{\Pi}_v\otimes\omega_{\pi_v}$. To do this we use some  
well-known facts regarding representations of $W_v$ and local Langlands  
correspondence for $\gl{n}(\rl)$ and $\gl{n}(\cx)$. We refer to  
\cite{knapp-motives} for a nice survey of these results.  
 
First assume that $k_v=\cx$. Then $W_v = \cx^\times$ and any irreducible  
representation of $W_v$ is one-dimensional and of the form  
\[ z \mapsto [z]^\ell |z|_\cx^t, \quad \ell\in\zl, t\in\cx, \]  
where $[z]=z/|z|$ and $|z|_\cx = |z|^2$.  
 
The $2n$-dimensional representation $\Phi_v$ of $W_v$ can now be written  
as a direct sum of $2n$ one-dimensional representations as above. Moreover,  
$\Phi_v(z) = \phi_v(z)$, considered as a diagonal matrix in $\gl{2n}(\cx)$,  
actually lies, up to conjugation, in $\gsp{2n}(\cx)$ or $\gso{2n}(\cx)$  
as in (\ref{torus-gsp-gso}).  
Therefore, there exist one-dimensional representations  
$\phi_0,\phi_1,\dots,\phi_n$ as above such that $\Phi_v$ is the direct  
sum of $\phi_1,\dots,\phi_n,\phi_n^{-1}\phi_0,\dots,\phi_1^{-1}\phi_0$. Now  
the central characters of $\Pi_v$ and $\pi_v$ can be written as  
$\oo_{\Pi_v}(z)=\det(\Phi_v(z))$ and $\oo_{\pi_v}(z) = e_0(\phi_v(z))$  
where $\phi_v(z)=\Phi_v(z)$ is considered as a $2n\times 2n$ diagonal  
matrix as in (\ref{torus-gsp-gso}) and $e_0$ is as in (\ref{e-e*}). In  
other words, $\oo_{\pi_v}=\phi_0$ and $\oo_{\Pi_v}=\phi_0^n$ or  
$\oo_{\Pi_v}=\oo_{\pi_v}^n$.  
 
Moreover, $\w{\Pi}_v$ corresponds to the $2n$-dimensional representation  
of $W_v$ which is the direct sum of $\phi_1^{-1},\dots,\phi_n^{-1}, 
\phi_n\phi_0^{-1},\dots,\phi_1\phi_0^{-1}$ implying that the two  
representations $\Pi_v$ and $\w{\Pi}_v\otimes\oo_{\pi_v}$ have the same  
parameters, i.e., $\Pi_v\simeq\w{\Pi}_v\otimes\oo_{\pi_v}$. 
 
Next assume that $k_v=\rl$. Then $W_v=\cx^\times \cup j\cx^\times$ with  
$j^2=-1$ and $j z j^{-1} = \bar{z}$ for $z\in\cx^\times$. Here the situation  
is identical and the only difference is that $W_v$ also has two-dimensional  
irreducible representations. The one-dimensional representations of $W_v$ can  
be described as  
\begin{eqnarray*}  
z \mapsto |z|_\rl^t,\quad j \mapsto 1, \quad t\in\cx \\ 
z \mapsto |z|_\rl^t,\quad j \mapsto -1, \quad t\in\cx  
\end{eqnarray*}  
with $|z|_\rl=|z|$ and the irreducible two-dimensional representations  
are of the form  
\[ z=r e^{i\theta} \mapsto  
\begin{matrix} 
\left( 
\begin{array}{cc}  
r^{2t} e^{i\ell\theta} & \\ 
 & r^{2t} e^{- i\ell\theta} 
\end{array} 
\right) 
\end{matrix}, \quad  
j \mapsto  
\begin{matrix} 
\left( 
\begin{array}{cc}  
0 & (-1)^\ell \\ 
1 & 0 
\end{array} 
\right) 
\end{matrix}, 
\]  
where $t\in\cx$ and $\ell\ge 1$ is an integer. These correspond, respectively,  
to representations $1\otimes|\cdot|_\rl^t$ and $\mbox{sgn}\otimes|\cdot|_\rl^t$  
of $\gl{1}(\rl)$ and $D_\ell\otimes|\cdot|_\rl^t$ of $\gl{2}(\rl)$ with  
notation as in \cite{knapp-motives}. 
 
Notice that again $\Phi_v(z)=\phi_v(z)$ is a diagonal $2n\times 2n$ matrix  
in $\gsp{2n}(\cx)$ or $\gso{2n}(\cx)$ as in the previous case while  
$\Phi_v(j)$ may have $2\times 2$ blocks as well. Therefore, we still have  
$\oo_{\Pi_v}=\oo_{\pi_v}^n$ and $\Pi_v \simeq\w{\Pi}_v\otimes\oo_{\pi_v}$  
as above.  
 
\item[{\bf(iii)}] $v<\infty$ and $\pi_v$ ramified:  
Choose $\Pi_v$ to be an arbitrary irreducible admissible   
representation of $\gl{2n}(k_v)$ with $\oo_{\Pi_v}=\oo_{\pi_v}^n$.  
\end{itemize} 
 
Then $\Pi = \otimes^\prime_v\  \Pi_v$ is an irreducible  
admissible representation of $\gl{2n}(\Ad)$ whose central character $\oo_\Pi$  
is equal to  $\oo_\pi^n$, and hence invariant under $k^\times$. Moreover,  
for all $v\not\in S$, we have that $L(s,\pi_v) = L(s,\Pi_v)$ by construction.  
Hence, $L^S(s,\Pi) =L^S(s,\pi)$, where  
\[ L^S(s,\Pi) = \prod_{v\not\in S} L(s,\Pi_v) \quad L^S(s,\pi) =  
\prod_{v\not\in S} L(s,\pi_v). \]  
Therefore, $L(s,\Pi) = \prod_v L(s,\Pi_v)$ is absolutely convergent in some  
right half plane. 
 
Choose $\eta=\otimes_v\ \eta_v$ to be a unitary character of  
$k^\times\backslash\Ad^\times$ such that $\eta_v$ is sufficiently ramified  
for $v\in S$ in order for Theorem \ref{stable} to hold and such that at one  
place $\eta_v^2$ is still ramified. For $\tau\in\cT(S;\eta)$ we claim the  
following equalities (along with their analogous equalities for the contragredients): 
\begin{eqnarray} 
  \label{l-equal}	L(s,\Pi_v\times\tau_v) &=& L(s,\pi_v\times\tau_v), \\ 
  \label{ep-equal}	\epsilon(s,\Pi_v\times\tau_v,\psi_v) &=&  
                        \epsilon(s,\pi_v\times\tau_v,\psi_v). 
\end{eqnarray} 
Here the $L$- and $\epsilon$-factors on the left are as in \cite{converse1}  
and those on the right are defined via the Langlands-Shahidi method  
\cite{shahidi:90annals,shahidi:88annals}. 
 
To see (\ref{l-equal}) and (\ref{ep-equal}) we again consider different  
places separately. 
 
\begin{itemize} 
 
\item[{\bf(i)}] $v<\infty$ and $\pi_v$ unramified:  
Let $\pi_v$ be again as in (\ref{UnramifiedChar}) with Satake parameter  
(\ref{parameter}). Then $\Pi_v$ will be as in (\ref{RepPi_v}). 
 
By \cite{jac-ps-shalika} we have  
\begin{eqnarray}\label{GLfactor} 
L(s,\Pi_v\times\tau_v) &=& \prod_{i=1}^n L(s,\tau_v\otimes\chi_i) 
L(s,\tau_v\otimes\chi_0\chi_i^{-1}) \\ 
\nonumber 
L(s,\w{\Pi}_v\times\w{\tau}_v) &=& \prod_{i=1}^n L(s,\w{\tau}_v\otimes\chi_i^{-1}) 
L(s,\w{\tau}_v\otimes\chi_0^{-1}\chi_i), 
\end{eqnarray} 
and  
\begin{equation}\label{GLfactor-ep} 
\e(s,\Pi_v\times\tau_v,\psi_v)=\prod_{i=1}^n \e(s,\tau_v\otimes\chi_i,\psi_v) 
                                  \e(s,\tau_v\otimes\chi_0\chi_i^{-1},\psi_v). 
\end{equation} 
 
On the other hand, it follows from the inductive property of $\g$-factors in  
the Langlands-Shahidi method (Theorem 3.5 of \cite{shahidi:90annals} or  
\cite{sha-psVolume}) that   
\begin{equation}\label{factor-gamma} 
\g(s,\pi_v\times\tau_v,\psi_v)=\prod_{i=1}^n \g(s,\tau_v\otimes\chi_i,\psi_v) 
                                  \g(s,\tau_v\otimes\chi_0\chi_i^{-1},\psi_v) 
\end{equation} 
just as in (\ref{mult}).

Since $\tau_v$ is generic, it is a full induced representation from generic  
tempered ones. Thus we can write  
\begin{equation} 
	\tau_v\simeq\mbox{Ind}(\nu^{b_1}\tau_{1,v}\otimes\cdots 
\otimes\nu^{b_p}\tau_{p,v}),  
\end{equation}  
where each $\tau_{j,v}$ is a tempered representation of some $\gl{r_j}(k_v)$,  
$\nu(\ )=\big|\det(\ )\big|_v$ on $\gl{r_j}(k_v)$, $r_1+\cdots + r_p = r$,  
and the $\tau_{j,v}$ are in the Langlands order. Moreover, recall that  
$\pi_v$ is the unique irreducible unramified  
subrepresentation of the representation of $\G_n(k_v)$ induced from the  
character $\chi$ as in (\ref{UnramifiedChar}) after an appropriate reordering,  
if necessary.  
 
Now by the definition of $L$-functions (Section 7 of \cite{shahidi:90annals})  
and their multiplicative property (Theorem 5.2 of \cite{sha-psVolume}) we have  
\begin{eqnarray}\label{factor} 
\nonumber 
  L(s,\pi_v\times\tau_v)&=&\prod_{j=1}^p L(s+b_j, \pi_v\times\tau_{j,v}) \\ 
\nonumber &=&\prod_{j=1}^p \prod_{i=1}^n L(s+b_j,\tau_{j,v}\otimes\chi_i)  
L(s+b_j,\tau_{j,v}\otimes\chi_0\chi_i^{-1})\\ 
          &=&\prod_{i=1}^n L(s,\tau_v\otimes\chi_i) L(s,\tau_v\otimes\chi_0\chi_i^{-1}). 
\end{eqnarray}				 
and likewise 
\begin{equation}\label{factor-contr} 
  L(s,\w{\pi}_v\times\w{\tau}_v)=\prod_{i=1}^n L(s,\w{\tau}_v\otimes\chi_i^{-1})  
L(s,\w{\tau}_v\otimes\chi_0^{-1}\chi_i). 
\end{equation} 
Note that Conjecture 5.1 of \cite{sha-psVolume}, which is a hypothesis of  
Theorem 5.2, is known in our cases by Theorem 5.7 of \cite{asgari}. 
 
Equations (\ref{factor-gamma}), (\ref{factor}), and (\ref{factor-contr}) in  
turn imply 
\begin{equation}\label{factor-ep} 
\e(s,\pi_v\times\tau_v,\psi_v)=\prod_{i=1}^n \e(s,\tau_v\otimes\chi_i,\psi_v) 
					     \e(s,\tau_v\otimes\chi_0\chi_i^{-1},\psi_v). 
\end{equation} 
 
Note that the product $L$-functions for $\gl{a}\times\gl{b}$ of the  
Langlands-Shahidi method  and the $L$-functions of \cite{jac-ps-shalika} are  
known to be equal (\cite{shahidi:84americanj}).  
Hence, to prove (\ref{l-equal}) and (\ref{ep-equal}) all we need is to  
compare the right hand sides of (\ref{GLfactor}) and (\ref{GLfactor-ep})  
with those of (\ref{factor}), (\ref{factor-contr}), and (\ref{factor-ep}). 
 
\item[{\bf(ii)}] $v|\infty$: By the local Langlands correspondence  
\cite{langlands-real} the representations $\pi_v$ and $\tau_v$ are given by  
admissible homomorphisms  
\[ \phi\,:\,W_v\longrightarrow\begin{cases} 
\gsp{2n}(\cx)\quad\mbox{ if } \G_n=\gspin{2n+1} \\  
\gso{2n}(\cx)\quad\mbox{ if } \G_n=\gspin{2n} \end{cases}  \] 
and  
\[ \phi^\prime\,:\,W_v\longrightarrow\gl{r}(\cx) \] 
respectively and the tensor product  
\[ (\iota\circ\phi)\otimes\phi^\prime\,:\,W_v\longrightarrow\gl{2nr}(\cx) \]  
is again admissible. Now,  
\[ L(s,\Pi_v\times\tau_v) = L(s,(\iota\circ\phi)\otimes\phi^\prime) =  
L(s,\pi_v\times\tau_v), \] 
and  
\[ \epsilon(s,\Pi_v\times\tau_v,\psi_v)=\epsilon(s,(\iota\circ\phi) 
\otimes\phi^\prime,\psi_v)=\epsilon(s,\pi_v\times\tau_v,\psi_v), \] 
where the middle factors are the local Artin-Weil factors \cite{tate-corvallis}  
and equalities hold by \cite{shahidi:85duke}. (See also \cite{borel-corvallis}.)  
 
\item[{\bf(iii)}] $v<\infty$ and $\pi_v$ ramified. This is where we will  
need the stability of $\g$-factors. Since $v\in S$ the representation  
$\tau_v$ can be written as  
\begin{equation} 
 \tau_v \simeq \mbox{Ind} (\nu^{b_1}\otimes\cdots\otimes\nu^{b_r})\otimes\eta_v  
        \simeq \mbox{Ind} (\eta_v\nu^{b_1}\otimes\cdots\otimes\eta_v\nu^{b_r}),  
\end{equation} 
where $\nu(x) = |x|_v$. Then,  
\begin{eqnarray} 
		\label{gspin-eqL} 
L(s,\pi_v\times\tau_v) &=& \prod_{i=1}^r L(s+b_i,\pi_v\times\eta_v), \\ 
\e(s,\pi_v\times\tau_v,\psi_v) &=& \prod_{i=1}^r \e(s+b_i,\pi_v\times\eta_v,\psi_v). 
\end{eqnarray} 
However, since $\eta_v$ is sufficiently ramified (depending on $\pi_v$),  
Corollary \ref{StableFormL-Ep} implies that  
\begin{eqnarray} 
L(s,\pi_v\times\eta_v)&\equiv&1 \\ 
\e(s,\pi_v\times\eta_v)&=&\prod_{i=1}^n \e(s,\eta_v\chi_i,\psi_v) 
\e(s,\eta_v\chi_0\mu_i^{-1},\psi_v)  
\end{eqnarray} 
for $n+1$ arbitrary characters $\chi_0,\chi_1,\dots,\chi_n$. We choose them  
to be as in (\ref{UnramifiedChar}). 
 
On the other hand, by either \cite{jac-ps-shalika} or \cite{shahidi:90annals} we have  
\begin{eqnarray} 
L(s,\Pi_v\times\tau_v) &=& \prod_{i=1}^r L(s+b_i,\Pi_v\otimes\eta_v), \\ 
\e(s,\Pi_v\times\tau_v,\psi_v) &=& \prod_{i=1}^r \e(s+b_i,\Pi_v\otimes\eta_v,\psi_v). 
\end{eqnarray} 
Again since $\eta_v$ is highly ramified (depending on $\Pi_v$) and  
$\oo_{\Pi_v}=\oo_{\pi_v}^n=\chi_0^n$ is equal to the  
product of the $2n$ characters  
\[ \chi_1,\dots,\chi_n,\chi_0\chi_n^{-1},\cdots,\chi_0\chi_1^{-1}, \]  
Proposition 2.2 of \cite{jac-shalika} implies that  
\begin{eqnarray} 
L(s,\Pi_v\otimes\eta_v)&\equiv&1 \\ 
\e(s,\Pi_v\otimes\eta_v)&=&\prod_{i=1}^n \e(s,\eta_v\chi_i,\psi_v) 
\e(s,\eta_v\chi_0\chi_i^{-1},\psi_v). 
			\label{gl-eqEP} 
\end{eqnarray} 
Comparing equations (\ref{gspin-eqL}) through (\ref{gl-eqEP}) now proves  
(\ref{l-equal}) and (\ref{ep-equal}) for $v$ non-archimedean with $\pi_v$ ramified. 
\end{itemize}

Now that we have (\ref{l-equal}) and (\ref{ep-equal}) for all places $v$ of $k$,  
we conclude globally that 
\begin{eqnarray} 
 \label{global-l} L(s,\Pi\times\tau) = L(s,\pi\times\tau) \quad&\quad 
                   L(s,\w{\Pi}\times\w{\tau}) = L(s,\w{\pi}\times\w{\tau}) \\ 
 \label{global-ep} \epsilon(s,\Pi\times\tau) = \epsilon(s,\pi\times\tau)  
  \quad&\quad \epsilon(s,\w{\Pi}\times\w{\tau}) = \epsilon(s,\w{\pi}\times\w{\tau})  
\end{eqnarray} 
for all $\tau\in\cT(S;\eta)$. All that remains now is to verify the three  
conditions of Theorem \ref{conversethm} which we can now check for the factors  
coming from the Langlands-Shahidi method thanks to (\ref{global-l}) and  
(\ref{global-ep}).  Conditions (1) -- (3) of Theorem \ref{conversethm} are  
Propositions \ref{entire}, \ref{bdd}, and \ref{fneq}, respectively. 
 
Therefore, there exists an automorphic representation $\Pi^\prime$ of  
$\gl{2n}(\Ad)$ such that for all $v\not\in S$ we have $\Pi_v\simeq\Pi^\prime_v$.  
In particular, for all $v\not\in S$ the local representation $\Pi^\prime_v$  
is related to $\pi_v$ as prescribed in Theorem \ref{transfer}.  
Moreover, note that for all $v\not\in S$ we have $\oo_{\Pi^\prime_v} =  
\oo_{\Pi_v} = \oo_{\pi_v}^n$. Since $\oo_{\Pi^\prime}$ is a gr{\"o}ssencharacter  
which agrees with the gr{\"o}ssencharacter $\oo_{\pi}^n$ at all but possibly  
finitely many places, we conclude that $\oo_{\Pi^\prime} = \oo_{\pi}^n$. 
 
On the other hand, if $v$ is an archimedean place or a non-archimedean place with  
$v\not\in S$, then we proved earlier that  
\[ \Pi^\prime_v \simeq \Pi_v \simeq \w{\Pi}_v\otimes\omega_{\pi_v}  
\simeq \w{\Pi}_v^\prime\otimes\omega_{\pi_v} \]  
which means, in particular, that $\Pi^\prime$ is nearly equivalent to  
$\w{\Pi}^\prime \otimes \omega_{\pi^\prime}$. 
\end{proof} 
 
\section{Complements}\label{apps} 
 
\subsection{Local Consequences} 
 
Our first local result is to show that the local transfers at  
the unramified places remain generic. Let us first recall a general  
result of Jian-Shu Li which we will use. The following is  
a special case of Theorem 2.2 of \cite{jshu}. 
 
\begin{prop}\label{local-generic}(J.-S. Li)  
Let $\G$ be a split connected reductive  
group over a non-archimedean local field $F$ and let $\B=\T\U$  
be a fixed Borel where $\T$ is a maximal torus and $\U$ is the  
unipotent radical of $\B$. Let $\chi$ be an unramified character  
of $\T(F)$ and let $\pi(\chi)$ be the unique irreducible unramified  
subquotient of the corresponding principal series representation.  
Then $\pi(\chi)$ is generic if and only if for all roots $\aa$  
of $(\G,\T)$ we have $\chi(\aa^\vee(\varpi))\not = |\varpi|_F$.  
Here $\aa^\vee$ denotes the coroot associated to $\aa$ and  
$\varpi$ is a uniformizer of $F$.  
\end{prop}

\begin{prop}\label{local-xfer} 
Let $\pi=\otimes_{v}^\prime\pi_v$ be an irreducible globally  
generic cuspidal automorphic representation of $\gspin{m}(\Ad)$,  
$m=2n+1$ or $2n$, and let $\Pi=\otimes_v^\prime\Pi_v$ be a  
transfer of $\pi$ to $\gl{2n}(\Ad)$ (cf. Theorem \ref{transfer}).  
If $v<\infty$ is a place of $k$ with $\pi_v$ unramified, then  
the local representation $\Pi_v$ is irreducible, unramified, and  
we have $\Pi_v\simeq\w{\Pi}_v\otimes\oo_{\pi_v}$. Moreover, if  
$m=2n+1$ (cf. Remark \ref{even-center}), then $\Pi_v$ is generic  
(and, hence, full induced principal series).  
\end{prop} 
 
\begin{proof} 
The representation $\Pi_v$ is irreducible and unramified by  
construction (cf. {\bf(i)} in the proof of Theorem \ref{transfer}).  
We also proved that $\Pi_v$ satisfies  
$\Pi_v\simeq\w{\Pi}_v\otimes\oo_{\pi_v}$ in the course of proof of  
Theorem \ref{transfer} in Section \ref{proof}. We now show that  
$\Pi_v$ is generic. Our tool will be Proposition \ref{local-generic}  
above.   
 
Now assume that $m=2n+1$.  
Let $\chi$ and $\chi_0, \dots, \chi_n$ be as in (\ref{UnramifiedChar}).  
Since $\pi_v$ is generic by Proposition \ref{local-generic} we have  
that $\chi(\aa^\vee(\varpi))\not = |\varpi|_{k_v}$ for all roots $\aa$.  
Using the notation of Section \ref{str}, the roots in the odd case $m=2n+1$  
are $\aa=\pm(e_i-e_j)$, $\pm(e_i+e_j)$ with $1\leq i \le j \leq n$ and  
$\pm(2 e_i)$ with $1\leq i \leq n$. The corresponding coroots are  
$\aa^\vee=\pm(e_i^*-e_j^*)$, $\pm(e_i^*+e_j^*-e_0^*)$ with  
$1\leq i \le j \leq n$ and $\pm(2 e_i^*-e_0^*)$ with $1\leq i \leq n$,  
respectively. This implies that  
$\chi_i \chi_j^{-1} \not= |\ |^{\pm 1}$ for $i\not= j$ and  
$\chi_i \chi_j \chi_0^{-1} \not= |\ |^{\pm 1}$ for all $i,j$. 
 
The representation $\Pi_v$ was chosen to be the unique irreducible  
unramified subquotient of the the representation on $\gl{2n}(F)$  
induced from the $2n$ unramified characters  
$\chi_1,\dots,\chi_n, \chi_0\chi_n^{-1},\dots,\chi_0\chi_1^{-1}$ as  
in (\ref{RepPi_v}). Therefore, the above relations imply that $\Pi_v$  
is generic and full-induced. 
\end{proof} 
 
\begin{rem}\label{even-center}  
The above argument does not quite work in the even case and one can  
easily construct local examples where the transferred local representation  
is the (unique) unramified subquotient of an induced representation  
on $\gl{2n}$ far from the generic constituent.  
 
For example, consider  
$\gspin{6}$ with $\chi_0=\mu^2$, $\chi_1 = \mu(5/2)$, $\chi_2=\mu(1/5)$,  
and $\chi_3=\mu(-3/2)$, where $\mu$ is a unitary character of $F^\times$  
and $\mu(r)$ means $\mu |\ |^r$. Now $\Pi_v$ is the unique unramified  
constituent of the representation on $\gl{6}(F)$ induced from  
$\mu(5/2)$, $\mu(3/2)$, $\mu(1/2)$, $\mu(-1/2)$, $\mu(-3/2)$, $\mu(-5/2)$ and, in fact,  
far from being generic. In this case, there is another constituent that  
is square-integrable, hence tempered and generic.  
 
Of course, we do expect $\Pi_v$ in the case of $m=2n$ to be generic as well.  
However, this phenomenon is not a purely local one in the case of $m=2n$.  
In fact, it will be automatic that the local transfers at  
the unramified places are generic once we prove that the automorphic representation  
$\Pi$ is induced from {\it unitary} cuspidal representations (see Remark  
\ref{unitary-remark} below). As we discuss in Remark \ref{unitary-remark} this  
will follow from our future work.  
\end{rem}

\subsection{Global Consequences}\label{global-sbsec} 
 
In this section we will make some comments about the automorphic  
representation $\Pi$ which are almost immediate consequences of our  
main result and leave more detailed information about $\Pi$ for a future  
paper.  
 
\begin{prop}\label{permute} 
Let $\pi$ be a globally generic cuspidal automorphic representation of  
$\gspin{m}(\Ad)$, $m=2n+1$ or $2n$ and let $\oo=\oo_\pi$. Then there  
exists a partition $(n_1,n_2,\dots,n_t)$ of $2n$ and (not necessarily unitary)  
cuspidal automorphic  
representations $\s_1,\dots,\s_t$ of $\gl{n_i}(\Ad)$, $i=1,\dots, t$,  
and permutation $p$ of $\{1,\dots,t\}$ with $n_i = n_{p(i)}$ and  
$\s_i\simeq \w{\s}_{p(i)} \otimes \oo$ such that any transfer $\Pi$ of $\pi$  
as in Theorem \ref{transfer} is a constituent of  
$\Sigma = \mbox{Ind} (\s_1\otimes\dots\otimes\s_t)$ where  
the induction is, as usual, from the standard parabolic of $\gl{2n}$ having  
Levi subgroup $\gl{n_1}\times\dots\times\gl{n_t}$.  
\end{prop} 
 
\begin{proof} 
Let $\Pi$ be any transfer of globally generic cuspidal representation  
$\pi$ as in Theorem \ref{transfer}. By Proposition 2 of \cite{langlands-notion}  
there exists a partition $p$ and $\s_i$'s as above such that $\Pi$ is a  
constituent of $\Sigma$. Furthermore, for finite places $v$ where $\pi_v$  
is unramified, we have that $\Pi_v$ is the unique unramified constituent  
of $\Sigma_v = \mbox{Ind} (\s_{1,v}\otimes\dots\otimes\s_{t,v})$. As part of  
Theorem \ref{transfer} we showed that $\Pi$ and $\w{\Pi}\otimes\oo$ are  
nearly equivalent (see definition prior to Theorem \ref{transfer}). Now,  
$\w{\Pi}\otimes\oo$ is a constituent of  
$\w{\Sigma}\otimes\oo=\mbox{Ind}((\w{\s}_1\otimes\oo)\otimes\cdots\otimes 
(\w{\s}_t\otimes\oo))$  
and by the classification theorem of Jacquet and Shalika (Theorem 4.4 of  
\cite{jac-sha-classification}) we have that there is a permutation $p$  
of $\{ 1,\dots,t\}$ such that $n_i = n_{p(i)}$ and  
$\s_{i}\simeq \w{\s}_{p(i)} \otimes \oo$.  
 
Now let $\Pi^\prime$ be another transfer of $\pi$ as in Theorem \ref{transfer}.  
Then, $\Pi^\prime$ is again a constituent of some  
$\Sigma^\prime = \mbox{Ind} (\s_{1}^\prime\otimes\dots\otimes\s_{t^\prime}^\prime)$,  
where $\s_{i}^\prime$  
is a cuspidal automorphic representation of $\gl{n_i^\prime}(\Ad)$  
and $(n_1^\prime,\dots,n^\prime_{t^\prime})$ is a partition of $2n$. Moreover,  
for almost all finite places $v$ we have that $\Pi^\prime_v$ is the unique  
unramified constituent of $\Sigma^\prime_v$. On the other hand, by construction,  
$\Pi_v \simeq \Pi^\prime_v$ for almost all $v$ and, therefore, the classification  
theorem of Jacquet and Shalika again implies that $t=t^\prime$ and, up to a  
permutation, $n_i=n^\prime_i$ and $\s_i\simeq \s^\prime_i$ for $i=1,\dots,t$.  
Therefore, $\Pi^\prime$ is also a constituent of $\Sigma$. 
\end{proof}  
 
\begin{rem}\label{unitary-remark}  
If we write $\s_i = \tau_i \otimes |\det(\ )|^{r_i}$ for $i=1,2,\dots,t$, with  
$\tau_i$ unitary cuspidal and $r_i\in\rl$, then we expect that all $r_i=0$, i.e.,  
$\Pi$ is an isobaric sum of unitary cuspidal representations. We will take up  
this issue, which will have important consequences, in our future work. 
\end{rem} 
 
\subsection{Exterior square transfer} 
 
In this section we show that exterior square generic transfer from  
$\gl{4}$ to $\gl{6}$ due to H. Kim (\cite{kim-ex2}) can be deduced  
as a special case of our main result. However, note that in this  
article we are only proving the weak transfer. Once we prove the  
strong version of transfer from $\gspin{2n}$ to $\gl{2n}$ again  
it will have the full content of the results of \cite{kim-ex2}. A  
similar remark also applies to Section \ref{gsp4}. 
 
\begin{prop}  
Let $\phi\,:\,\gspin{6}\vlto{30}\gl{4}$ be the (double) covering map  
(cf. Proposition \ref{spin-structure}) and denote by $\h{\phi}$ the map  
induced on the connected components of the $L$-groups. 
\[ \begin{array}{ccl} 
\gso{6}(\cx) & \overset{\iota}{\vlto{20}} & \gl{6}(\cx) \\ 
&& \\ 
{\h{\phi}\nto{20}} & \neto{25} & \\ 
\gl{4}(\cx) & & 
\end{array} \] 
 
Then, the $\iota\circ\h{\phi} = \wedge^2$. 
\end{prop}  
\begin{proof} 
The group $\gso{6}$ is of type $D_3$ and we denote its simple roots by  
$\a{1},\a{2},\a{3}$ as in Section \ref{str}. Also, $\gl{4}$ is of  
type $A_3$ (or $D_3$) and we denote its corresponding simple roots by  
$\ba{2},\ba{1},\ba{3}$, respectively, and similarly for other root data  
(cf. Section \ref{str}).  
Let $A=\text{diag}(a_1,a_2,a_3,a_4)\in\gl{4}(\cx)$.  
For a fixed appropriate choice of fourth root of unity and  
$\delta=(a_1 a_2 a_3 a_4)^{1/4}$ we have  
\begin{eqnarray*} 
\iota\circ\h{\phi}(A) &=& \iota\circ\h{\phi}\,  
(\, \delta \, \ba{2}^\vee(\frac{a_1}\delta) \ba{1}^\vee(\frac{a_1a_2}{\delta^2})  
\ba{3}^\vee(\frac{a_1 a_2 a_3}{\delta^3}) \, ) \\ 
&=& \iota(e_0^*(\delta^4)e_1^*(\delta^2)e_2^*(\delta^2)e_3^*(\delta^2)  
\a{2}^\vee(\frac{a_1}\delta) \a{1}^\vee(\frac{a_1a_2}{\delta^2})  
\a{3}^\vee(\frac{a_1 a_2 a_3}{\delta^3})) \\ 
&=& \iota(e_0^*(\delta^4) e_1^*(a_1a_2) e_2^*(a_1a_3) e_3^*(a_2a_3)) \\ 
&=& \text{diag}(a_1a_2, a_1a_3,a_2a_3, a_2a_4, a_1a_4,a_3a_4)=\wedge^2 A.  
\end{eqnarray*} 
Here the third equality follows from Proposition \ref{cent-gs}. 
\end{proof} 
 
As a corollary we see that our Theorem \ref{transfer} in the special case of  
$m=2n$ with $n=3$ gives Kim's exterior square transfer. 
 
\begin{prop} If $\pi$ is an irreducible cuspidal automorphic  
representation of $\gl{4}(\Ad)$ considered as a representation of  
$\gspin{6}(\Ad)$ via the covering map $\phi$, then the automorphic  
representation $\Pi$ of Theorem \ref{transfer} is such that  
$\Pi_v = \wedge^2 \pi_v$ for almost all $v$. 
\end{prop} 
 
\subsection{Transfer from $\gsp{4}$ to $\gl{4}$}\label{gsp4} 
 
The special case of $m=2n+1$ with $n=4$ of our Theorem \ref{transfer}  
gives the following: 
 
\begin{prop} 
Let $\pi$ be an irreducible globally generic cuspidal automorphic  
representation of $\gsp{4}(\Ad)$. Then $\pi$ can be transferred  
to an automorphic representation $\Pi$ of $\gl{4}(\Ad)$ associated  
to the embedding $\gsp{4}(\cx)\hookrightarrow\gl{4}(\cx)$.  
\end{prop} 
\begin{proof} 
Notice that $\gspin{5}$ is isomorphic, as an algebraic group, to  
the group $\gsp{4}$. Now the corollary is a special case of Theorem \ref{transfer}  
as mentioned above. 
\end{proof} 

\begin{rem}\label{ramin}
The above Proposition, in particular, proves that the spinor $L$-function of 
$\pi$ is entire. We understand that R. Takloo-Bighash also has 
a proof of this fact using a completely different method based on integral 
representations. 
\end{rem} 
 
\bibliography{gspin} 
\bibliographystyle{plain}

\end{document}